\documentclass[journal]{IEEEtran}

\ifCLASSINFOpdf

\else

\fi

\usepackage[cmex10]{amsmath}
\usepackage{graphics,graphicx,epic,eepic,epsfig,times,units} %
\usepackage{psfrag}
\usepackage{resizegather}		
\usepackage{cite}				
\usepackage{xcolor}				
\usepackage{cases}				
\usepackage{subcaption}			

\usepackage{alphalph}
\usepackage{etoolbox}

\AtBeginDocument{%
	\AtBeginEnvironment{subequations}{%
		\let\alph\alphalphval%
	}
}

\usepackage[hyphens]{url}
\usepackage[hidelinks]{hyperref}			
\usepackage{breakurl}
\hypersetup{breaklinks=true}

\usepackage{cleveref}			
\crefformat{equation}{(#2#1#3)}
\crefrangeformat{equation}{(#3#1#4)--(#5#2#6)}
\crefmultiformat{equation}{(#2#1#3)}%
{,~(#2#1#3)}{, (#2#1#3)}{,~(#2#1#3)}

\DeclareMathOperator*{\argmin}{argmin}

\usepackage{tabularx}
\usepackage{multirow}
\usepackage{xcolor}
\usepackage{siunitx}
\DeclareSIUnit \VAr {VAr} 
\DeclareSIUnit \VA {VA} 
\DeclareSIUnit \rad {Radians} 

\newcommand{\reals}{{\mbox{\bf R}}}

\newcommand{\dom}{{\mbox{\bf dom}}}
\newcommand{\epi}{{\mbox{\bf epi}}}

\usepackage{algorithm,algorithmicx}
\usepackage{algpseudocode}
\usepackage{algcompatible}
\algrenewcommand\alglinenumber[1]{\scriptsize #1:}
\newcommand{\algrule}[1][.2pt]{\par\vskip.5\baselineskip\hrule height #1\par\vskip.5\baselineskip}
\algrenewcommand\algorithmicindent{2.0em}%

\usepackage{changepage}

\usepackage{mathtools,nccmath,bm}

\usepackage[font=small,skip=0pt]{caption}

\usepackage{amsthm,amssymb}
\newtheorem{theorem}{Theorem}

\newtheorem{definition}{Definition}
\let\olddefinition\definition
\renewcommand{\definition}{\olddefinition\normalfont}

\allowdisplaybreaks

\begin{document}

\title{An Exact Sequential Linear Programming Algorithm for the Optimal Power Flow Problem}

\author{\IEEEauthorblockN{Sleiman~Mhanna,~\emph{Member, IEEE},
		Pierluigi Mancarella,~\emph{Senior~Member, IEEE}}}

\maketitle

\begin{abstract}
	Despite major advancements in nonlinear programming (NLP) and convex relaxations, most system operators around the world still predominantly use some form of linear programming (LP) approximation of the AC power flow equations. This is largely due to LP technology's superior reliability and computational efficiency, especially in real-time market applications, security-constrained applications, and extensions involving integer variables, in addition to its ability to readily generate locational marginal prices (LMP) for market applications. In the aim of leveraging the advantages of LP while retaining the accuracy of NLP interior-point methods (IPMs), this paper proposes a sequential linear programming (SLP) approach consisting of a sequence of carefully constructed \emph{supporting hyperplanes and halfspaces}. The algorithm is numerically demonstrated to converge on 138 test cases with up the 3375 buses to feasible high-quality solutions (i) \emph{without} AC feasibility restoration (\emph{i.e.,} using LP solvers exclusively), (ii) in computation times generally within the same order of magnitude as those from a state-of-the-art NLP solver, and (iii) with robustness against the choice of starting point. In particular, the (relative) optimality gaps and the mean constraint violations are on average around $10^{-3}$\% and $10^{-7}$, respectively, under a \emph{single} parameter setting for all the 138 test cases. To the best of our knowledge, the proposed SLP approach is the first to use LP exclusively to reach \emph{feasible} and high-quality solutions to the nonconvex AC OPF in a reliable way, which paves the way for system and market operators to keep using their LP solvers but now with the ability to accurately capture transmission losses, price reactive power (Q-LMP), and obtain more accurate LMP.
\end{abstract}

\IEEEpeerreviewmaketitle

\setlength{\belowdisplayskip}{0.5pt} \setlength{\belowdisplayshortskip}{0.5pt}
\setlength{\abovedisplayskip}{0.5pt} \setlength{\abovedisplayshortskip}{0.5pt}

\section*{Notation}
\addcontentsline{toc}{section}{Notation}
\subsection{Sets}
	\begin{IEEEdescription}[\IEEEsetlabelwidth{$\Re\left\{\bullet\right\}$}\IEEEusemathlabelsep]
		\item[$\mathcal{B}$] Set of all buses in the power network.
		\item[$\mathcal{B}_{i}$] Set of buses connected to bus $i$.
		\item[$\mathcal{G}$] Set of all generators $gi$  such that $i$ is the bus and $g$ is the generator connected to it.
		\item[$\mathcal{G}_{i}$] Set of all generators connected to bus $i$.
		\item[$\mathcal{L}$] Set of all branches $ij$ where $i$ is the ``from'' bus.
		\item[$\mathcal{L}_{t}$] Set of all branches $ji$ where $j$ is the ``to'' bus.
	\end{IEEEdescription}
\subsection{Parameters and input data}
	\begin{IEEEdescription}[\IEEEsetlabelwidth{$\Re\left\{\bullet\right\}$}\IEEEusemathlabelsep]
		\item[$b^{\rm sh}_{i}$] Shunt susceptance (pu) at bus $i$.
		\item[$b^{\rm ch}_{ij}$] Charging susceptance (pu) in the $\pi$-model of line $ij$.
		\item[$c_{0,gi}$] Constant coefficient ($\SI{}{\$\per\hour}$) term of generator $g$'s cost function.
		\item[$c_{1,gi}$] Coefficient ($\SI{}{\$\per\mega\watt\hour}$) of the linear term of generator $g$'s cost function.
		\item[$c_{2,gi}$] Coefficient ($\SI{}{\$\per\mega\watt\hour\squared}$) of the quadratic term of generator $g$'s cost function.
		\item[$g^{\rm sh}_{i}$] Shunt conductance (pu) at bus $i$.
		\item[$k$] Iteration number.
		\item[$p_{i}^{\rm d}/q_{i}^{\rm d}$] Active/Reactive power demand (pu) at bus $i$.
		\item[$\overline{s}_{ij}$] Apparent power rating (pu) of branch $ij$.
		\item[$\overline{\theta}_{ij}/\underline{\theta}_{ij}$] Upper/Lower limit of the difference of voltage angles of buses $i$ and $j$.
		\item[$\rho_{ij}$] Penalty parameter used to guide convergence.
		\item[$\tilde{t}_{ij}$] Complex tap ratio of a phase shifting transformer ($\tilde{t}_{ij}=\tau_{ij}\mathrm{e}^{\mathrm{i} \theta_{i}^{\rm shift}}$).
		\item[$\tilde{y}_{ij}$] Series admittance (pu) in the $\pi$-model of line $ij$.
		\item[$\zeta$] Percentage apparent power loading of a branch.
	\end{IEEEdescription}
\subsection{Operators}
	\begin{IEEEdescription}[\IEEEsetlabelwidth{$\Re\left\{\bullet\right\}$}\IEEEusemathlabelsep]
		\item[$\bullet^*$] Conjugate operator.
		\item[$\dom f$] Domain of function $f$.
		\item[$\epi f$] Epigraph of function $f$.
		\item[$\nabla f$] Gradient of function $f$.
		\item[$\Im/\Re\left\{\bullet\right\}$] Imaginary/Real value operator.
		\item[$\underline{\bullet}/\overline{\bullet}$] Minimum/Maximum magnitude operator.
		\item[$\left|\bullet\right|$] Magnitude operator/Cardinality of a set.
		\item[$\bullet^T$] Transpose operator.
	\end{IEEEdescription}
\subsection{Variables}
	\begin{IEEEdescription}[\IEEEsetlabelwidth{$\Re\left\{\bullet\right\}$}\IEEEusemathlabelsep]
		\item[$p_{gi}/q_{gi}$] Active/Reactive power (pu) generation of generator $g$ at bus $i$. 
		\item[$p_{ij}/q_{ij}$] Active/Reactive power (pu) flow on branch $ij$.
		\item[$r_{ij}$] Non-negative slack variable used to (i) guide convergence and (ii) circumvent infeasbility.
		\item[$\theta_{i}$] Voltage angle ($\SI{}{\radian}$) at bus $i$.
		\item[$v_{i}$] Voltage magnitude (pu) at bus $i$.
	\end{IEEEdescription}
\subsection{Acronyms}
	\begin{IEEEdescription}[\IEEEsetlabelwidth{$\Re\left\{\bullet\right\}$}\IEEEusemathlabelsep]
		\item[AC] Alternating current. 
		\item[AEMO] Australian Energy Market Operator.
		\item[DC] Direct current.
		\item[IPM] Interior-point method.
		\item[LMP] Locational marginal prices/pricing.
		\item[LP] Linear programming.
		\item[MILP] Mixed-integer linear programming.
		\item[NLP] Nonlinear programming.
		\item[NZEM] New Zealand's Electricity Market.
		\item[OPF] Optimal power flow.
		\item[QCP] Quadratically constrained programming.
		\item[Q-LMP] Locational marginal reactive power prices/pricing.
		\item[RTO] Regional Transmission Organizations.
		\item[SDP] Semidefinite programming.	
		\item[SLP] Sequential linear programming.
		\item[SOCP] Second-order cone programming.	
	\end{IEEEdescription}

\section{Introduction}

	\IEEEPARstart{T}{he} classical optimal power flow (OPF) problem consists of finding the least-cost dispatch of active and reactive power from generators to satisfy the loads at all buses in a way that is governed by physical laws, such as Ohm's law and Kirchhoff's law, as well as other technical restrictions such as thermal limit constraints. The main complexity of the OPF problem is directly attributed to definition of complex power in which the product of the complex voltage and the conjugate of the complex current gives rise to the nonlinear and nonconvex alternating current (AC) power flow constraints. Consequently, finding a globally optimal solution to this nonconvex nonlinear programming (NLP) problem is proven to be non-deterministic polynomial-time hard (NP-hard) \cite{Lehmann2016_ACfeasibility,Bienstock2015_NPhardnessofACPF}. Since its inception by Carpentier in 1962 \cite{Carpentier1962_OPF}, the OPF problem has garnered rigorous research attention that gave rise to a rich body of books and literature on solution approaches that can be broadly classified into (i) exact methods, (ii) approximations, and (iii) convex relaxations.
	
	Exact methods solve the AC OPF directly to find a \emph{feasible}, globally or locally optimal solution. The set of exact methods guaranteeing locally optimal solutions includes primal-dual IPMs with a trust-region and/or filter line search \cite{Momoh1999_ReviewofOPFII,Momoh1999_ReviewofOPFI,Frank2012_OPFSurveypartI,Capitanescu2016_Criticalreview,MATPOWER}, sequential quadratic programming \cite{KNITRO}, and sequential linear programming \cite{Castillo2016_SLPforOPF}. The computationally more demanding global solution methods include \emph{single-tree} methods such as spatial branch and bound \cite{Gopalakrishnan2012_BandBfortheOPF}, and \emph{multi-tree} methods such as outer approximation \cite{Liu2018_MultitreeforGlobalOPF}. State-of-the-art IPMs are computationally efficient and have been shown to solve large-scale instances (more than 50,000 buses) in minutes \cite{Kardos2018_NumericalEvalofIPMforOPF}. Nonetheless, despite the tractability of today's IPMs, most system and market operators around the world still favor LP approximations such as the direct current (DC) OPF for the following reasons: a) IPMs are generally sensitive to the choice of starting point \cite{Stott2012_OPFbasicrequirements} and to the choice of formulation \cite{Coffrin2020_ImpactofPWLcostonOPF} (especially in cases involving piecewise linear generator cost functions), which means that the computation time can go from seconds to minutes, b) LP solvers are reliable and always converge to a \emph{unique} solution, c) Locational Marginal Prices (LMP), which are the linchpin of electricity markets, are easier to obtain using LP, d) The computational superiority of LP and mixed-integer LP (MILP) solvers is particularly advantageous in computationally demanding areas such as real-time markets, security-constrained OPF, day-ahead security-constrained unit commitment, and transmission network expansion planning to name a few \cite{Stott2009_DCPFrevisited}, e) In complex engineering systems such as the power network, where controls are constantly changing, the concept of optimality (local or global) may be elusive, if not completely irrelevant, compared to other more basic requirements \cite{Stott2012_OPFbasicrequirements}. In other words, the industry traded exactness for reliability, robustness, and computational efficiency.
	
	Because the DC OPF inherently ignores losses, a common technique in academic literature is to incorporate losses by approximating the cosine terms in the original polar-form AC power flow constraints by their Maclaurin series expansion and then linearizing \cite{DosSantos2011_PWLforDCOPFlosses} or convexifing \cite{Zhong2013_EDwithQuadraticLosses} the resulting nonconvex quadratic constraints. On the other hand, system operators around the world adopt different approaches to capture active power (MW) losses. The Australian Energy Market Operator (AEMO) and New Zealand's Electricity Market (NZEM) model the losses as \emph{marginal loss factors} in the form of quadratic functions of real power, which are then approximated by piecewise linear segments \cite{Hobbs2008_QuadraticLossOPF}. Some system operators and Regional Transmission Organizations (RTO) in the USA use \emph{loss factors} together with \emph{loss distribution factors} to represent the sensitivity of the losses with respect to power injections at each bus \cite{Yang2018_ImprovingMWonlyOPF}. System operators in China distribute an estimate of the losses to the loads \cite{Yang2018_ImprovingMWonlyOPF}. 
	
	Unfortunately, even after modeling the losses, the DC OPF can result in a poor approximation of the true solution and the associated LMP \cite{Stott2009_DCPFrevisited}. In fact, it is mathematically proven in \cite{Baker2019_DCOPFisneverfeasible} that solutions from the DC OPF are never AC feasible. In an attempt to restore AC feasibility, some system operators in the USA resort to subsequent corrective measures in the form of a quasi-optimization process that iterates between a DC OPF with loss factors and an AC power flow \cite{ONeill2011_FERConISOsoftware}. Furthermore, because the DC OPF ignores reactive power, system operators and RTO is the USA do not price reactive power, which leads to \emph{uplift}, whereby out-of-market costs are incurred by out-of-merit resources to relieve network constraints and provide reactive power services \cite{FederalEnergyRegulatoryCommission2014_Uplift}.
	
	Another research strand revolves around convex relaxations which generally aim to find the \emph{convex hull} of the nonconvex feasible region in the hope that the solution of the resulting convex problem is exact, in which case the global optimum of an NP-hard problem can be computed in polynomial time. Convex relaxations trace back to \cite{Jabr2006_SOCforOPF}, which numerically showed that a second-order cone programming (SOCP) relaxation of the \emph{alternative-form} OPF \cite{Esposito1999_Reliableloadflow} \emph{can} be exact on radial networks, and to \cite{Bai2008_SDPforOPF}, which numerically showed that a semidefinite programming (SDP) relaxation \emph{can} be exact on meshed networks. Since then, it has been well understood that convex relaxations are only exact under specific conditions that are unfortunately difficult to predict in practice and can only be truly verified \emph{after} the problem solved \cite{Low2014_ConvexOPF_Exactness,Huang2017_SufficientConditioninDistributionNetworks}. This has therefore inspired a myriad of solution approaches on tightening techniques towards recovering feasible solutions to the original AC OPF with guaranteed global optimality in the best case, and computing tighter lower bounds in the worst case (see \cite{Zohrizadeh2019_SOCP_OPF} for a comprehensive review). Alas, despite the reliability of convex programming technology, convex relaxations still have not found their way into real-world industrial applications as they generally lack the computational efficiency of LP-based approximations.
	
	In light of the above, one way to capture both losses and reactive power while still using LP technology is by using a carefully designed SLP approach. SLP, which was originally introduced in \cite{Cheney1959_NewtonsmethodandTchebycheffApprox,Kelley1960_CuttingPlaneMethod} as the \emph{cutting plane method}, consists of solving the original NLP problem by solving a series of LP problems generated by approximating all the nonlinear constraints by their first-order Taylor series expansion. Recently, \cite{Castillo2016_SLPforOPF} proposed an SLP approach using the \emph{IV formulation} and numerically demonstrated its potential in solving OPF instances with up to 3375 buses to high-quality solutions. The method controls the step-size by introducing a penalty function and corresponding slack variables. However, rigorous numerical analysis in \cite{Castillo2016_SLPforOPF} revealed the sensitivity of the method to the choice of starting point and the authors indeed acknowledge this being one of the pitfalls of the approach. Another major limitation of applying an SLP scheme to the IV formulation is the eventuality of dealing with nonconvex \emph{inequality} constraints (the nonconvex quadratic voltage bounds) whose linearization cuts the interior of the feasible region, which increases the risk of infeasible LP problems and suboptimal solutions. In an attempt to circumvent the shortcomings of \cite{Castillo2016_SLPforOPF}, the work in \cite{Sampath2018_SLPforOPF} uses a trust-region based SLP approach on the \emph{polar-form} OPF but now with an AC feasibility restoration phase consisting of a Newton-Raphson scheme to bring the iterates close to the boundary of the feasible region. However, the work in \cite{Sampath2018_SLPforOPF} requires problem-specific parameter tuning and the numerical evaluation shows that the method converges to poor-quality solutions with a marginal improvement in computation time. The approach in \cite{Yang2016_SLPforOPF} also uses a second-order AC feasibility restoration phase but now with a combination of careful mathematical transformations and first-order approximations of the bilinear voltage terms instead of a step-size control. A related work by the same authors proposes an improved DC OPF model that captures both reactive power and voltages (and therefore losses) in the aim of bestowing reactive power pricing to the incumbent AC-DC quasi-optimization processes adopted by most system operators and RTO in the USA \cite{Yang2018_ImprovingMWonlyOPF}.

	Against this background, this work attempts to further narrow the gap between LP approximations and IPMs by proposing an exact SLP method that leverages the structure of the \emph{alternative form} \cite{Esposito1999_Reliableloadflow} of the OPF problem by dynamically introducing affine and linear cuts as a sequence of carefully constructed supporting hyperplanes and an increasing number supporting halfspaces. The success of this particular construction hinges on a \emph{single} parameter, which is automatically tuned. This is in contrast to the trust-region-like approaches which may require introducing a lot more slack variables and parameters that require cumbersome and sometimes problem-specific tuning. Another distinctive feature of the proposed method is that it does not require an AC feasibility restoration phase.
	
	The method is numerically demonstrated to converge to feasible high-quality solutions in computation times within the same order of magnitude as a state-of-the-art IPM solver, without intermediate AC feasibility restoration and with robustness of the choice of starting point, under a single parameter setting for all the test cases. In particular, the (relative) optimality gaps and the mean constraint violations are on average around $10^{-3}$\% and $10^{-7}$, respectively. Furthermore, the method is not only evaluated on meshed network instances from MATPOWER \cite{MATPOWER} (as is the case in \cite{Castillo2016_SLPforOPF}, \cite{Yang2016_SLPforOPF}, and \cite{Sampath2018_SLPforOPF}), but also on the more challenging ones from PGLib-OPF \cite{Babaeinejadsarookolaee2019_PGlibOPF} as well as on radial distribution network instances, for a total of 138 test cases. Finally, because it uses only LP solvers to directly solve the AC OPF to feasibility, the proposed SLP method can also be viewed as a promising alternative to the widely used quasi-optimization process that alternates between DC OPF with loss factors and AC feasibility, which can lead to uplift, costlier operation, and inaccurate LMP.
	
	In a nutshell, this paper advances the state of the art in the following ways:
	\begin{itemize}
		\item It introduces a novel SLP algorithm for the nonconvex AC OPF problem consisting of a sequence of carefully constructed supporting hyperplanes and halfspaces. The exactness, computational efficiency, and robustness against the choice of starting point are demonstrated in an extensive numerical analysis on 138 test cases with both radial and meshed topologies with up to 3375 buses.
		\item In addition to LMP, the method can readily generate accurate Q-LMP using LP solvers exclusively, which lays the foundation for reactive power pricing.
	\end{itemize}
	
	The paper is organized as follows. Section \ref{Sec_OPF} describes the alternative form of the OPF problem as well as the SOCP relaxation. Section \ref{sec_SLPforOPF} introduces three novel SLP approaches in Sections \ref{subsec_SLPforSOCR}, \ref{subsec_SLPforOPFRadial}, and \ref{subsec_SLPforOPFMeshed} for the SOCP relaxation of the OPF, the OPF in radial networks, and the OPF in general (meshed) networks, respectively. Section \ref{sec_numericalevalulation} numerically evaluates the proposed algorithms for general networks on 138 instances with up to 3375 buses. The paper concludes in Section \ref{sec_conclusion}. 

\section{Optimal power flow}\label{Sec_OPF}

	Consider a power network $\mathcal{N}=\left(\mathcal{B},\mathcal{L} \cup \mathcal{L}_{\rm t} \right) $, where $\mathcal{B}$ is the set of buses and $\mathcal{L} \cup \mathcal{L}_{\rm t}$ is the set of edges such that $\mathcal{L}$ is the set of all branches $ij$ where $i$ is the ``from'' bus, and $\mathcal{L}_{\rm t}$ is the set of all branches $ji$ where $j$ is the ``to'' bus. Examples of power system branches include transmission lines and phase-shifting transformers. There are two types of power network topologies, \emph{radial} and \emph{meshed}. Medium and low-voltage distribution networks are generally radial, whereas high-voltage transmission networks are almost exclusively meshed. The set of generating units is denoted by the two-dimensional tuple set $\mathcal{G}=\left\lbrace \left( g,i \right) | g \in \mathcal{G}_{i} , i \in \mathcal{B} \right\rbrace $, where $\mathcal{G}_{i}$ is the set of all generators connected to bus $i$. In the classical AC OPF, the demand at corresponding buses is assumed to be \emph{static}. The complex voltage $\tilde{v}_{i}$ (pu) at bus $i$ can be expressed as $\tilde{v}_{i}=v_{i} \mathrm{e}^{\mathrm{i} \theta_{i}}=v_{i} \angle \theta_{i} = v_{i}\cos\left(\theta_{i}\right) + \mathrm{i} v_{i}\sin\left(\theta_{i}\right)$ in \emph{polar form}, where $\mathrm{i} = \sqrt{-1}$. In power system analysis, the voltage magnitude is usually nondimensionalized against a base voltage level and is expressed in per-unit (pu). For instance, if the base voltage is $\SI{100}{\kilo\volt}$, then the pu equivalent of $\SI{110}{\kilo\volt}$ would be 1.1 (pu). In practice, bus voltage magnitudes $v_{i}$ at the transmission system level are desired to be close to the base voltage to ensure system security and stability. Transmission lines and phase-shifting transformers are represented by their $\pi$-model equivalents, in which the admittance is defined by $\tilde{y}_{ij}=g_{ij} + \mathrm{i} b_{ij}$, where $g_{ij}$ and $b_{ij}$ are the conductance (pu) and susceptance (pu), respectively. Additionally, the charging susceptance in the $\pi$-model of branch $ij$ is denoted by $b^{\rm ch}_{ij}$ (pu).
	
	By defining
	\begin{subequations}\label{opf_definingW}
		\begin{align}
			w_{i} & =\left|\tilde{v}_{i}\right|^2 = v^{2}_{i}, & i \in \mathcal{B}, \\
			w_{ij}^{\rm r} & =\Re\left\{\tilde{v}_{i} \tilde{v}_{j}^*\right\}=v_{i}v_{j}\cos\left(\theta_{i}-\theta_{j}\right), & ij \in \mathcal{L}, \\
			w_{ij}^{\rm i} & =\Im\left\{\tilde{v}_{i} \tilde{v}_{j}^*\right\}=v_{i}v_{j}\sin\left(\theta_{i}-\theta_{j}\right), & ij \in \mathcal{L},
		\end{align}
	\end{subequations}
	the \emph{alternative formulation} \cite{Esposito1999_Reliableloadflow} of the OPF problem can be written as
	\begin{subequations}\label{opf_alt}	
		\begin{align}
			& \hspace{-1.5cm} \underset {\substack{p_{gi},q_{gi},w_{i},\theta_{i},\\w_{ij}^{\rm r},w_{ij}^{\rm i},p_{ij},q_{ij}}} 
			{\mbox{ minimize}} \quad \sum_{gi \in \mathcal{G}} f_{gi}\left(p_{gi}\right) & & \label{opf_alt_objective} \\
			& \hspace{-1.5cm} \text{ subject to} & & \nonumber \\
			\underline{p}_{gi} & \leq p_{gi} \leq \overline{p}_{gi}, & \hspace{-2.5cm} gi\in \mathcal{G} & \label{opf_alt_Pminmax} \\
			\underline{q}_{gi} & \leq q_{gi} \leq \overline{q}_{gi}, & \hspace{-2.5cm} gi\in \mathcal{G} & \label{opf_alt_Qminmax} \\
			\underline{v}_{i}^2 & \leq w_{i} \leq \overline{v}_{i}^2, & i \in \mathcal{B} & \label{opf_alt_Vminmax} \\ 
			\underline{\theta}_{ij} & \leq \theta_{i}-\theta_{j} \leq \overline{\theta}_{ij}, & ij \in \mathcal{L} & \label{opf_alt_Angle} \\
			\sum_{g \in \mathcal{G}_{i}} p_{gi}&-p_{i}^{\rm d} = \sum_{j \in \mathcal{B}_{i}} p_{ij} + g^{\rm sh}_{i}w_{i}, & i \in \mathcal{B} & \label{opf_alt_kclp} \\	
			\sum_{g \in \mathcal{G}_{i}} q_{gi}&-q_{i}^{\rm d} = \sum_{j \in \mathcal{B}_{i}} q_{ij} - b^{\rm sh}_{i}w_{i}, & i \in \mathcal{B} & \label{opf_alt_kclq} \\
			p_{ij}=& \ g^{\rm c}_{ij} w_{i} - g_{ij}w_{ij}^{\rm r} + b_{ij}w_{ij}^{\rm i}, & ij \in \mathcal{L} \cup \mathcal{L}&_{\rm t} & \label{opf_alt_pij} \\
			q_{ij}=& \ b^{\rm c}_{ij} w_{i} -b_{ij} w_{ij}^{\rm r} - g_{ij}w_{ij}^{\rm i}, & ij \in \mathcal{L} \cup \mathcal{L}&_{\rm t} & \label{opf_alt_qij} \\
			w_{i}w_{j} & = (w_{ij}^{\rm r})^2 + (w_{ij}^{\rm i})^2, & ij \in \mathcal{L} & \label{opf_alt_RSOC} \\
			\theta_{i} -\theta_{j} & = {\rm atan2}(w_{ij}^{\rm i},w_{ij}^{\rm r}), & ij \in \mathcal{L} & \label{opf_alt_atan2}\\
			& p_{ij}^2+q_{ij}^2 \leq \overline{s}_{ij}^2, & ij \in \mathcal{L} \cup \mathcal{L}&_{\rm t} \label{opf_alt_thermal}
		\end{align}
	\end{subequations}
	where $g^{\rm c}_{ij}=\Re\{(\tilde{y}_{ij}^*-0.5\mathrm{i}b^{\rm ch}_{ij})/\left|\tilde{t}_{ij}\right|^2\}$, $b^{\rm c}_{ij}=\Im\{(\tilde{y}_{ij}^*-0.5\mathrm{i}b^{\rm ch}_{ij})$ $ /\left|\tilde{t}_{ij}\right|^2\}$, $g_{ij}=\Re\{\tilde{y}_{ij}^*/\tilde{t}_{ij}\}$, $b_{ij}=\Im\{\tilde{y}_{ij}^*/\tilde{t}_{ij}\}$, $g^{\rm c}_{ji}=\Re\{\tilde{y}_{ji}^*-0.5\mathrm{i}b^{\rm ch}_{ji}\}$, $b^{\rm c}_{ji}=\Im\{\tilde{y}_{ji}^*-0.5\mathrm{i}b^{\rm ch}_{ji}\}$, $g_{ji}=\Re\{\tilde{y}_{ji}^*/\tilde{t}_{ji}^*\}$, $b_{ji}=\Im\{\tilde{y}_{ji}^*/\tilde{t}_{ji}^*\}$, and $\tilde{t}_{ij}=\tau_{ij}\mathrm{e}^{\mathrm{i} \theta_{ij}^{\rm shift}}$ is the complex tap ratio of a phase shifting transformer.\footnote{For a transmission line $\tau_{ij}=1$ and $\theta_{ij}^{\rm shift}=0$.} The objective function in \eqref{opf_alt_objective} is generally a quadratic cost function of the form $f_{gi}\left(p_{gi}\right)=c_{2,gi}\left(p_{gi}\right)^2+c_{1,gi}\left(p_{gi}\right)+c_{0,gi}$. The relationship between $w_{i}$, $w_{ij}^{\rm r}$ and $w_{ij}^{\rm i}$ in \eqref{opf_definingW} is captured by \cref{opf_alt_RSOC}, which is a rotated second-order cone (SOC) constraint first introduced in \cite{Esposito1999_Reliableloadflow}.\footnote{Note that $w_{ji}^{\rm i} = -w_{ij}^{\rm i}$ and $w_{ji}^{\rm r} = w_{ij}^{\rm r}$.} The nonconvexities in \eqref{opf_alt} stem from equality constraint \cref{opf_alt_RSOC}, which describes the \emph{boundary} of a rotated SOC, and \cref{opf_alt_atan2} which contains the nonconvex ${\rm atan2}$ function. In meshed networks, \cref{opf_alt_atan2} is necessary to ensure that angle differences sum to zero over \emph{every} cycle in the network. In radial networks, which are acyclic, constraint \cref{opf_alt_atan2} can be dropped from \eqref{opf_alt} and the angles can be recovered by solving a system of linear equations of the form \cref{opf_alt_atan2} from the optimal $w_{ij}^{\rm r}$ and $w_{ij}^{\rm i}$. Specifically, the OPF in balanced 3-phase radial distribution networks can be written as
	\begin{subequations}\label{opf_radial}	
		\begin{align}
			& \hspace{-1.25cm} \underset {\substack{p_{gi},q_{gi},w_{i},\\w_{ij}^{\rm r},w_{ij}^{\rm i},p_{ij},q_{ij}}} 
			{\mbox{ minimize}} \quad \sum_{gi \in \mathcal{G}} f_{gi}\left(p_{gi}\right) & & \label{opf_radial_objective} \\
			& \hspace{-1.25cm} \text{ subject to \cref{opf_alt_Pminmax,opf_alt_Qminmax,opf_alt_Vminmax}, \cref{opf_alt_kclp,opf_alt_kclq,opf_alt_pij,opf_alt_qij,opf_alt_RSOC},\cref{opf_alt_thermal}}. & & \label{opf_radial_shared}
		\end{align}
	\end{subequations}
	The OPF formulation in \eqref{opf_alt} is equivalent to the polar, rectangular, and IV OPF formulations (see \cite{ONeill2012_IVformulation} for a review). 
	
	The alternative formulation of the OPF problem possesses an appealing property. In particular, because voltages are non-negative (\emph{i.e.,} $w_{i} \ge 0 \ \forall i \in \mathcal{B}$), relaxing \cref{opf_alt_RSOC} into an inequality constraint makes it convex. In fact, this property was the intuition behind the well-known SOCP relaxation in \cite{Jabr2006_SOCforOPF}, which further ignores \cref{opf_alt_atan2}. As such, the SOCP relaxation of the OPF can be written as
	\begin{subequations}\label{opf_socr}	
		\begin{align}
			& \hspace{-1.25cm} \underset {\substack{p_{gi},q_{gi},w_{i},\\w_{ij}^{\rm r},w_{ij}^{\rm i},p_{ij},q_{ij}}} 
			{\mbox{ minimize}} \quad \sum_{gi \in \mathcal{G}} f_{gi}\left(p_{gi}\right) & & \label{opf_socr_objective} \\
			& \hspace{-1.25cm} \text{ subject to \cref{opf_alt_Pminmax,opf_alt_Qminmax,opf_alt_Vminmax}, \cref{opf_alt_kclp,opf_alt_kclq,opf_alt_pij,opf_alt_qij},\cref{opf_alt_thermal}} & & \label{opf_socr_shared} \\
			w_{i}w_{j} & \ge (w_{ij}^{\rm r})^2 + (w_{ij}^{\rm i})^2, & ij \in \mathcal{L} & \label{opf_socr_RSOC} \\
			\tan{\underline{\theta}_{ij}}w_{ij}^{\rm r} & \le w_{ij}^{\rm i} \le \tan{\overline{\theta}_{ij}}w_{ij}^{\rm r}. & ij \in \mathcal{L} \label{opf_SOC_Angle}
		\end{align}
	\end{subequations} 
	Problem \eqref{opf_socr} can be written in the general form 
	\begin{subequations}\label{opf_socr_Generalproblem}
		\begin{align}
			\underset {\substack{x}} 
			{\mbox{ minimize}} \quad f_{0}\left(x\right) & & \label{opf_socr_Generalproblem_objective}\\
			\text{ subject to} \quad f_{i}\left(x\right) & \leq 0, & & i=1,\ldots,m \label{opf_socr_Generalproblem_eq_elec}\\
			a_{i}^{T}x & \leq b_{i}, & & i=1,\ldots,p \label{opf_socr_Generalproblem_linear}
		\end{align}
	\end{subequations}
	where $x \in \reals^{n}$, $f_{0}$,$f_{1},\ldots,f_{m}: \reals^{n} \rightarrow \reals$, $a_{1},\ldots,a_{p} \in \reals^{n}$, and $b_{1},\ldots,b_{p} \in \reals$. Function $f_{0}\left(x\right)$ represents the cost functions in \eqref{opf_socr_objective}, whereas \eqref{opf_socr_Generalproblem_eq_elec} represent the convex nonlinear constraints \eqref{opf_socr_RSOC} and \eqref{opf_alt_thermal}. The linear inequality constraints in \eqref{opf_socr_Generalproblem_linear} represent \cref{opf_alt_Pminmax,opf_alt_Qminmax,opf_alt_Vminmax} and \cref{opf_alt_kclp,opf_alt_kclq,opf_alt_pij,opf_alt_qij}.\footnote{Recall that an affine equality constraint of the form $a_{i}^{T}x = b_{i}$ can be rewritten as $\left\lbrace a_{i}^{T}x \leq b_{i}\right\rbrace \cap \left\lbrace - a_{i}^{T}x \leq -b_{i} \right\rbrace $.}
	
	In light of the above, this work capitalizes on this specific property of the alternative formulation to design an SLP approach incorporating a sequence of carefully constructed supporting hyperplanes \emph{in conjunction with} an increasing number of halfspaces bounded by the supporting hyperplanes of a nonconvex set of the form
	\begin{multline}\label{QoL}
		\mathcal{W}_{ij} = \left\lbrace \left(w_{i},w_{j},w_{ij}^{\rm r},w_{ij}^{\rm i}\right) \in \reals^{4} | \right. \\
		\left. w_{i} = \left( (w_{ij}^{\rm r})^2 +( w_{ij}^{\rm i})^2 \right) /w_{j}, w_{i},w_{j} \ge 0 \right\rbrace, ij \in \mathcal{L} .
	\end{multline}
	A similar treatment cannot be applied to \cref{opf_alt_atan2} as this constraint does not become convex when the equality is relaxed into an inequality. However, \cref{opf_alt_atan2} is quasilinear in $\theta_{i}-\theta_{j}$ and $w_{ij}^{\rm i}/w_{ij}^{\rm r}$,\footnote{This is because $ w_{ij}^{\rm r} \in [\underline{v}_{i}\underline{v}_{j}\cos{\underline{\theta}_{ij}}, \overline{v}_{i}\overline{v}_{j} ] $ and $ w_{ij}^{\rm i} \in [\overline{v}_{i}\overline{v}_{j}\sin{\underline{\theta}_{ij}}, \overline{v}_{i}\overline{v}_{j}\sin{\overline{\theta}_{ij}} ] $, and therefore the \emph{range} of ${\rm atan2}(w_{ij}^{\rm i},w_{ij}^{\rm r})$ lies in a quasilinear region.} and therefore a simple first-order approximation is enough to approximate it to a very high accuracy in an iterative LP-based algorithm, especially when used in tandem with the above construction for the set $\mathcal{W}_{ij}$. This elaborate construction, which is detailed in the next section, obviates the need for a step-size control in the classical sense (see \cite{Castillo2016_SLPforOPF} and \cite{Sampath2018_SLPforOPF}).
	
	Another advantage of the alternative formulation (Problem~\ref{opf_alt}) of the AC OPF problem is that it enjoys \emph{linear} branch current limit constraints. The derivation of those linear constraints can be found in Appendix~\ref{sec_currentlimit}. Some system operators deal directly with branch current limits, as opposed to branch MVA limits, since the current magnitude is directly related to conductor temperature \cite{Shchetinin2019_CurrentLimits}. Although, those linear branch current limit constraints can be directly handled by the proposed SLP approach without extra treatment, the branch MVA limits are used (and then linearized) instead for generality and comparison purposes.

\section{Sequential linear programming}\label{sec_SLPforOPF}

	This section introduces three novel SLP approaches in Sections \ref{subsec_SLPforSOCR}, \ref{subsec_SLPforOPFRadial}, and \ref{subsec_SLPforOPFMeshed} for the SOCP relaxation of the OPF (Problem \ref{opf_socr}), the nonconvex AC OPF in radial networks (Problem \ref{opf_radial}), and the nonconvex AC OPF in meshed networks (Problem \ref{opf_alt}), respectively. Instead of directly introducing the SLP algorithm for the full nonconvex AC OPF problem in \eqref{opf_alt}, which may be overwhelming at first glance, the below specific sequence of introductions aims at easing the exposition by walking the reader through the evolution of the SLP algorithm from one that solves the convex SOCP relaxation of the OPF problem in \eqref{opf_socr} to one that solves the original nonconvex AC OPF problem in \eqref{opf_alt}.
	
\subsection{SLP for the SOCP relaxation of the OPF}\label{subsec_SLPforSOCR}
	
	Before introducing an SLP approach for Problem \ref{opf_socr}, constraint \eqref{opf_socr_RSOC} is rewritten as 
	\begin{align}\label{opf_socr_RSOC_Slei}
		& w_{i} \ge \frac{(w_{ij}^{\rm r})^2+(w_{ij}^{\rm i})^2}{w_{j}}, & ij \in \mathcal{L} &
	\end{align}
	which is also a convex set as it can be interpreted as the \emph{epigraph} of a convex \emph{quadratic-over-linear} function $f_{ij}: \reals^{3} \rightarrow \reals$ defined as
	\begin{align*}
		\epi f_{ij} = \left\lbrace \left( w,w_{i} \right) | w \in \dom f_{ij}, f_{ij}(w) \le w_{i} \right\rbrace,
	\end{align*}
	where $f_{ij}(w) = \left( (w_{ij}^{\rm r})^2+(w_{ij}^{\rm i})^2\right) /w_{j}$, $\dom f_{ij} = \reals^{2} \times \reals_{++} $, and $w = [w_{ij}^{\rm r},w_{ij}^{\rm i},w_{j}]^{T}$. This form is more conducive to convergence as now the linearization is in the reduced space of 3 variables instead of 4 (when using \eqref{opf_socr_RSOC} directly). The first-order Taylor series approximation of $f_{ij}(w)$ at a point $w^{(k)}$ can be written as
	\begin{multline*}
		f_{ij}^{\rm a}(w,w^{(k)}) = f_{ij}\left( w^{(k)}\right) + \nabla f_{ij}\left( w^{(k)}\right)^{T}\left(w - w^{(k)} \right) \\
		= f_{ij}\left( w^{(k)}\right) + 
		\begin{bmatrix}
			2w_{ij}^{{\rm r},(k)}/w_{j}^{(k)} \\
			2w_{ij}^{{\rm i},(k)}/w_{j}^{(k)} \\
			\displaystyle \frac{\left(w_{ij}^{{\rm r},(k)}\right)^2+\left(w_{ij}^{{\rm i},(k)}\right)^2}{-\left(w_{j}^{(k)}\right)^2}
		\end{bmatrix}^{T}
		\begin{bmatrix}
			w_{ij}^{{\rm r}} - w_{ij}^{{\rm r},(k)} \\
			w_{ij}^{{\rm i}} - w_{ij}^{{\rm i},(k)} \\
			w_{j}-w_{j}^{(k)}
		\end{bmatrix}.
	\end{multline*}
	It is easy to show that the hyperplane $\{ (w,w_{i}) | w_{i} - f_{ij}^{\rm a}(w,w^{(k)}) = 0 \} $ \emph{supports} the set $\epi f_{ij}$ at the boundary point $\left(w^{(k)}, f_{ij}( w^{(k)} ) \right) $, since for any
	\begin{align*}
	\resizebox{\linewidth}{!}{$
		\left( w,w_{i} \right) \in \epi f_{ij} \Rightarrow
		\begin{bmatrix}
			\nabla f_{ij}\left( w^{(k)}\right) \\
			-1
		\end{bmatrix}^{T}
		\begin{pmatrix}
			\begin{bmatrix}
				w \\
				w_{i}
			\end{bmatrix}-
			\begin{bmatrix}
				w^{(k)} \\
				 f_{ij}\left( w^{(k)}\right)
			\end{bmatrix}
		\end{pmatrix}
			\le 0.$}
	 \end{align*}
	Said differently, the set $\{ (w,w_{i}) | w_{i} \ge f_{ij}^{\rm a}(w,w^{(k)}) \} $ is a \emph{supporting halfspace} to the set 
	$\{ (w,w_{i}) | w_{i} \ge \left( (w_{ij}^{\rm r})^2+(w_{ij}^{\rm i})^2\right) /w_{j}, w_{j} > 0 \} $ at the boundary point $\left(w^{(k)},\right. $ $ \left. f_{ij}( w^{(k)} ) \right) $. 
	
	The second step consists of approximating the convex quadratic thermal limit constraint in \cref{opf_alt_thermal}. Instead of directly approximating \cref{opf_alt_thermal} by a first-order Taylor series expansion around $( p_{ij}^{(k)},q_{ij}^{(k)}) $, the approximation is applied to the \emph{projection} of point $( p_{ij}^{(k)},q_{ij}^{(k)} ) $ on the (convex) set defined by \cref{opf_alt_thermal}. In more detail, let $\mathcal{C} = \{ y \in \reals^{2} | \left\| y \right\|^{2}_{2} \le \overline{s}^{2}_{ij} \}$, where $y = \left[ p_{ij},q_{ij}\right]^{T} $. The projection of $y^{(k)}$ on $\mathcal{C}$, denoted by $P_{\mathcal{C}}\left( y^{(k)} \right) $, can now be obtained by solving the following problem
	\begin{align*}
		P_{\mathcal{C}}\left( y^{(k)} \right) = \argmin_{y}\left\lbrace \left\| y - y^{(k)} \right\|^{2}_{2} | y \in \mathcal{C} \right\rbrace,
	\end{align*} 
	which admits the analytical solution
	\begin{align*}
		\resizebox{\linewidth}{!}{$P_{\mathcal{C}}\left( y^{(k)} \right) = \left[ \frac{\overline{s}_{ij}p_{ij}^{(k)} }{\sqrt{ \left(p_{ij}^{(k)}\right)^{2} + \left(q_{ij}^{(k)}\right)^{2} }}, \frac{\overline{s}_{ij}q_{ij}^{(k)} }{\sqrt{ \left(p_{ij}^{(k)}\right)^{2} + \left(q_{ij}^{(k)}\right)^{2} }} \right]^{T}.$} 
	\end{align*} 
	After factorizing and rearranging the terms, the first-order Taylor series approximation of $g_{ij}\left( y \right) = p_{ij}^{2} + q_{ij}^{2} - \overline{s}^{2}_{ij}$ at the boundary point $P_{\mathcal{C}}\left( y^{(k)} \right)$ can be written as
	\begin{align*}
		g_{ij}^{\rm a}\left( y,y^{(k)} \right) = -2\overline{s}^{2}_{ij} + 2y^{T}P_{\mathcal{C}}\left( y^{(k)} \right).
	\end{align*}
	The set $\{ y | g_{ij}^{\rm a}\left( y,y^{(k)} \right) \le 0 \} $ is a supporting halfspace to set $\mathcal{C}$ at the boundary point $P_{\mathcal{C}}\left( y^{(k)} \right) $. 
	
	The last step is to substitute the quadratic terms in $f_{gi}\left(p_{gi}\right)$ by corresponding variables and rotated SOC constraints which can then be tightly approximated by a lifted polyhedron as in \cite{Mhanna2016_TightLPfortheOPF}. This construction requires far less inequality constraints than standard piecewise linear (PWL) approximations for the same accuracy. Specifically, $\left|\mathcal{G}\right|$ variables $p_{g}$ and associated linear constraints of the form $p_{g}=\sqrt{c_{2,gi}}p_{gi}$, for all $ gi \in \mathcal{G}$ are introduced along with $N=\left\lfloor \left|\mathcal{G}\right| /2\right\rfloor+\left\lceil \left|\mathcal{G}\right| /2 - \left\lfloor \left|\mathcal{G}\right| /2\right\rfloor \right\rceil$ variables $\alpha_{n}$ and constraints of the form $\alpha_{n} \geq p_{2n-1}^2+p_{2n}^2$ for which the lifted polyhedral construction is denoted as
	\begin{align}\label{RSOC_Obj}
		\mathcal{P}_{k}^{\rm R} \left( p_{2n-1}^2+p_{2n}^2 \leq \alpha_{n} \right), \qquad n \in \left\{1,\ldots,N\right\}.
	\end{align}
	
	A polyhedral outer approximation of Problem \ref{opf_socr} can now be iteratively constructed by dynamically adding supporting halfspaces to the LP approximation of Problem \ref{opf_socr} at iteration $k-1$. The SLP method for solving the SOCP relaxation of the OPF in Problem \ref{opf_socr} is described in Algorithm~\ref{SLP_SOCP}, where
	\begin{align*}
		\hspace{-0.2cm} F = \left[ \left( w_{i}^{(k)}-\left( \left( w_{ij}^{{\rm r},(k)} \right)^2 + \left( w_{ij}^{{\rm i},(k)} \right) ^2\right) /w_{j}^{(k)} \right)_{ij \in \mathcal{L}} \right],
	\end{align*}
	and 
	\begin{align*}
		G = \left[ \left( \left( p_{ij}^{(k)} \right)^{2} + \left( q_{ij}^{(k)} \right)^{2} - \overline{s}^{2}_{ij} \right)_{ij \in \mathcal{L} \cup \mathcal{L}_{\rm t}} \right].
	\end{align*}
	\begin{algorithm}[!t]
		\caption{SLP for the SOCP relaxation of the OPF.}
		\begin{adjustwidth}{-1em}{}
			\begin{algorithmic}[1]
				\scriptsize
				\STATE \parbox[t]{\dimexpr\linewidth-0.75cm}{\textbf{Initialization:}
					Choose starting point $\left(v_{i}^{(0)},\theta_{i}^{(0)}\right)$ for all $i \in \mathcal{B}$. Set $k:=1$, $\epsilon,\epsilon^{\rm th} \in \left[ 10^{-10}, 10^{-3} \right] $, $\overline{k} = 50 $, $\zeta \in (0,1)$, and $\mathcal{T}_{ij} := \left\lbrace \emptyset \right\rbrace$ for all $ ij \in \mathcal{L} \cup \mathcal{L}_{\rm t} $. \strut}
				\algrule
				\WHILE {$\min \left( F \right) < -\epsilon$ \textbf{and} ${\rm max}\left( G \right) > \epsilon^{\rm th}$ \textbf{and} $ k \le \overline{k}$}
				\STATE \parbox[t]{\dimexpr\linewidth-0.75cm}{Obtain $\left( w_{i}^{(k)},w_{ij}^{{\rm r},(k)},w_{ij}^{{\rm i},(k)}, p_{ij}^{(k)}, q_{ij}^{(k)} \right) $ by solving
					\begin{subequations}\label{opf_socr_slp}	
						\begin{align}
							& \hspace{-0.75cm} \underset {\substack{p_{gi},p_{g},\alpha_{n},q_{gi},w_{i},\\w_{ij}^{\rm r},w_{ij}^{\rm i},p_{ij},q_{ij}}} 
							{\mbox{ minimize}} \quad \sum_{n=1}^{N} \alpha_{n} + \sum_{gi \in \mathcal{G}} c_{1,gi}\left(p_{gi}\right) + c_{0,gi} \hspace{-3cm} & & \label{opf_socr_slp_objective} \\
							& \hspace{-0.5cm} \text{ subject to \cref{opf_alt_Pminmax,opf_alt_Qminmax,opf_alt_Vminmax}, \cref{opf_alt_kclp,opf_alt_kclq,opf_alt_pij,opf_alt_qij}, \cref{opf_SOC_Angle}, \cref{RSOC_Obj}} \hspace{-2cm} & & \label{ opf_socr_slp_shared} \\
							& f_{ij}^{\rm a} \left( w,w^{(\kappa)} \right) \le w_{i}, & ij \in \mathcal{L}, \ \kappa=0,\ldots,k-1, & \label{opf_socr_slp_cut1} \\
							& g_{ij}^{\rm a}\left( y,y^{(\kappa)} \right) \le 0, \hspace{-1cm} & ij \in \mathcal{L} \cup \mathcal{L}_{\rm t}, \ \kappa \in \mathcal{T}_{ij}, & \label{opf_socr_slp_cut2}
						\end{align}
					\end{subequations} \strut} 
				\FOR {$ij \in \mathcal{L} \cup \mathcal{L}_{\rm t} $}
					\IF {$ \left( \left( p_{ij}^{(k)} \right)^{2} + \left( q_{ij}^{(k)} \right)^{2} \right) > \left( \zeta \overline{s}_{ij} \right)^{2} $}
						\STATE $\mathcal{T}_{ij} := \mathcal{T}_{ij} \cup \{k\}$.
					\ENDIF
				\ENDFOR
				\STATE {$k:= k + 1 $.}
				\ENDWHILE
			\end{algorithmic} 
		\end{adjustwidth}
		\label{SLP_SOCP}
	\end{algorithm}
	At every iteration $k$, Algorithm~\ref{SLP_SOCP} solves an LP approximation of Problem \ref{opf_socr} in the form of Problem \ref{opf_socr_slp} (line 3), which features at most $\left| \mathcal{L} \right| + \left| \mathcal{L} \cup \mathcal{L}_{\rm t} \right| $ additional linear constraints compared to the problem at $k-1$. These supporting halfspaces, delineated by \eqref{opf_socr_slp_cut1} and \eqref{opf_socr_slp_cut2}, are dynamically added to Problem \ref{opf_socr_slp} at $k$. Set $\mathcal{T}_{ij}$ registers all previous iteration numbers that qualify as halfspaces for \eqref{opf_alt_thermal}. The condition on line 5 is satisfied when constraint \eqref{opf_alt_thermal} is $(1-\zeta^2)\overline{s}_{ij}^2$ away from being violated, in which case the supporting halfspace $\{ y | g_{ij}^{\rm a}\left( y,y^{(k)} \right) \le 0 \} $ is added to Problem \ref{opf_socr_slp} at $k$. In other words, the condition on line 5 is true when the apparent power flow on branch $ij$ is greater than $\zeta\overline{s}_{ij}$, where $\zeta$ is a parameter that captures a predetermined percentage loading of a branch. The algorithm terminates when the maximum violations of constraints \eqref{opf_socr_RSOC_Slei} and \eqref{opf_alt_thermal} do not exceed the small tolerances $\epsilon$ and $\epsilon^{\rm th}$. Aside from these tolerances and $\zeta$, Algorithm~\ref{SLP_SOCP} requires no other parameters to ensure convergence. In fact, since the SOCP relaxation of the OPF is a convex problem, the SLP in Algorithm~\ref{SLP_SOCP} with an arbitrary starting point that lies in the domain of all the constraint functions is guaranteed to converge to an optimal solution of Problem \ref{opf_socr} in a finite number of iterations (if Problem \ref{opf_socr} is feasible and a solution exists).
	\begin{theorem}\label{opf_socr_slp_Convergence}
		Let the feasible set of Problem \ref{opf_socr_Generalproblem} be defined as $\mathcal{G} = \{ x | f_{i}\left(x\right) \leq 0, i=1,\ldots,m, a_{i}^{T}x \leq b_{i}, i=1,\ldots,p \} $, which is nonempty and compact. Problem \ref{opf_socr_Generalproblem} then consists of finding a vector $x^{\star}$ such that $f_{0}(x^{\star})=\min \left\lbrace f_{0}(x) | x \in \mathcal{G} \right\rbrace $. If $x^{(k)} \in \mathcal{S}^{(k)}$ is such that
		\begin{align}
			f_{0}(x^{(k)})=\min \left\lbrace f_{0}(x) | x \in \mathcal{S}^{(k)} \right\rbrace, 
		\end{align}
		where $\mathcal{S}^{(0)} = \left\lbrace x | a_{i}^{T}x \leq b_{i}, i=1,\ldots,p \right\rbrace $ is compact and 
		\begin{align*}
			\mathcal{S}^{(k)}= \mathcal{S}^{(k-1)} \cap \left\lbrace x | f_{i}(x^{(k)}) + \nabla f_{i}(x^{(k)})^{T}(x - x^{(k})\leq 0, \right. \\
			\left. i=1,\ldots,m \right\rbrace, 
		\end{align*}
		then the sequence $x^{(k)}$ contains a Cauchy subsequence that converges to a point $x^{\star} \in \mathcal{G}$ with $f_{0}(x^{(0)}) \le f_{0}(x^{(1)}) \le \cdots \le f_{0}(x^{(k-1)}) \le f_{0}(x^{(k)}) \le f_{0}(x^{\star}) \le f_{0}(x)$ for all $x \in \mathcal{G}$. 
	\end{theorem}
	The proof can be found in Appendix~\ref{opf_socr_slp_Proof}. Since Algorithm~\ref{SLP_SOCP} terminates when $\min \left( F \right) \ge -\epsilon$ and ${\rm max}\left( G \right) \le \epsilon^{\rm th}$, $x^{\star}$ is called an $\epsilon-$effective solution. Finally, although the proof of Theorem~\ref{opf_socr_slp_Convergence} depended on it, not all the supporting halfspaces in \eqref{opf_socr_slp_cut1} are necessary for convergence. Previous inactive supporting halfspaces can be carefully removed in exchange for a reduction in the overall computation time. 
	
	Although quadratically constrained programming (QCP) solvers such as Gurobi \cite{Gurobi2019} and CPLEX \cite{CPLEX} can efficiently handle large-scale problems of the form \eqref{opf_socr}, the \emph{Barrier method} is widely known to be numerically more stable on large-scale LP than on large-scale QCP problems \cite{Gurobi2019,Mhanna2016_TightLPfortheOPF}. More importantly, in applications requiring binary or integer decision variables, another advantage of using Algorithm~\ref{SLP_SOCP} is that a \emph{crossover} strategy can be used to recover an \emph{LP basis} that can be used in the branch-and-cut algorithms at the core of MILP solvers (as these use the primal or dual simplex which allow for efficient re-optimization at each node of the branch-and-bound tree), which are more efficient than both the branch-and-cut with the linearized outer approximation and the branch-and-bound with IPM (usually Barrier) techniques used in mixed-integer QCP (MIQCP) solvers. 
	
	Under certain conditions, the SOCP relaxation can be exact in radial networks \cite{Low2014_ConvexOPF_Exactness,Huang2017_SufficientConditioninDistributionNetworks}, \emph{i.e.,} constrains \eqref{opf_socr_RSOC_Slei} are all active at the optimum, in which case Algorithm~\ref{SLP_SOCP} can be directly used in a sequential MILP framework when binary or integer variables are present. Moreover, Algorithm~\ref{SLP_SOCP} can easily be adapted to the SOC relaxation of the DistFlow model \cite{Baran1989_DistFlow} for radial networks. In the next section, Algorithm~\ref{SLP_SOCP} is extended to a general setting where the SOCP relaxation may not be exact in 3-phase balanced radial distribution networks.
	
\subsection{Proposed SLP algorithm for the OPF in radial networks}\label{subsec_SLPforOPFRadial}
	
	\begin{algorithm}[!t]
		\caption{SLP for the OPF in radial networks.}
		\begin{adjustwidth}{-1em}{}
			\begin{algorithmic}[1]
				\scriptsize
				\STATE \parbox[t]{\dimexpr\linewidth-0.75cm}{\textbf{Initialization:}
					Choose starting point $\left(v_{i}^{(0)},\theta_{i}^{(0)}\right)$ for all $i \in \mathcal{B}$. Set $k:=1$, $\epsilon,\epsilon^{\rm th} \in \left[ 10^{-10}, 10^{-3} \right] $, $\overline{k} = 50 $, $\zeta \in (0,1)$, $ \overline{\rho} > \rho_{ij}^{(0)} > 0$, $\gamma > 0 $, $\mathcal{I}_{ij} := \left\lbrace \emptyset \right\rbrace$ for all $ ij \in \mathcal{L}$, and $\mathcal{T}_{ij} := \left\lbrace \emptyset \right\rbrace$ for all $ ij \in \mathcal{L} \cup \mathcal{L}_{\rm t} $. \strut}
				\algrule
				\WHILE {$\max \left( F \right) > \epsilon$ \textbf{and} ${\rm max}\left( G \right) > \epsilon^{\rm th}$ \textbf{and} $ k \le \overline{k}$}
				\STATE \parbox[t]{\dimexpr\linewidth-0.75cm}{Obtain $\left( w_{i}^{(k)},w_{ij}^{{\rm r},(k)},w_{ij}^{{\rm i},(k)}, p_{ij}^{(k)}, q_{ij}^{(k)} \right) $ by solving
					\begin{subequations}\label{SLP_radial}	
						\begin{align}
							& \hspace{-0.75cm} \underset {\substack{p_{gi},p_{g},\alpha_{n},q_{gi},w_{i},\\w_{ij}^{\rm r},w_{ij}^{\rm i},p_{ij},q_{ij},r_{ij}}} 
							{\mbox{ minimize}} \eqref{opf_socr_slp_objective} + \sum_{ij \in \mathcal{L}} \rho_{ij}^{(k-1)}r_{ij} \hspace{-3cm} & & \label{SLP_radial_objective} \\
							& \hspace{-0.5cm} \text{ subject to \cref{opf_alt_Pminmax,opf_alt_Qminmax,opf_alt_Vminmax}, \cref{opf_alt_kclp,opf_alt_kclq,opf_alt_pij,opf_alt_qij}, \cref{opf_SOC_Angle}, \cref{RSOC_Obj}} \hspace{-3cm} & & \label{SLP_radial_shared} \\
							& f_{ij}^{\rm a} \left( w,w^{(k-1)} \right) + r_{ij} = w_{i}, & ij \in \mathcal{L}, & \label{SLP_radial_cut0} \\
							& f_{ij}^{\rm a} \left( w,w^{(\kappa)} \right) \le w_{i}, & ij \in \mathcal{L}, \ \kappa \in \mathcal{I}_{ij}, & \label{SLP_radial_cut1} \\
							& g_{ij}^{\rm a}\left( y,y^{(\kappa)} \right) \le 0, \hspace{-2cm} & ij \in \mathcal{L} \cup \mathcal{L}_{\rm t}, \ \kappa \in \mathcal{T}_{ij}, & \label{SLP_radial_cut_thermal} \\
							& r_{ij} \ge 0, & ij \in \mathcal{L}, & \label{opf_slp_r}
						\end{align}
					\end{subequations} \strut} 
				\FOR {$ij \in \mathcal{L} $}
					\IF {$ F_{ij} \notin \left[-\epsilon, \epsilon \right] $}
						\STATE { $\mathcal{I}_{ij} := \mathcal{I}_{ij} \cup \{k-1\} $.}
					\ENDIF
					\IF {$ r_{ij} \ge \epsilon $}
						\STATE $\rho_{ij}^{(k)} := \min\left( \overline{\rho}, \gamma \rho_{ij}^{(k-1)} \right) $.
					\ENDIF
				\ENDFOR
				\FOR {$ij \in \mathcal{L} \cup \mathcal{L}_{\rm t} $}
					\IF {$ \left( \left( p_{ij}^{(k)} \right)^{2} + \left( q_{ij}^{(k)} \right)^{2} \right) > \left( \zeta \overline{s}_{ij} \right)^{2} $}
						\STATE $\mathcal{T}_{ij} := \mathcal{T}_{ij} \cup \{k\}$.
					\ENDIF
				\ENDFOR
				\STATE {$k:= k + 1 $.}
				\ENDWHILE
			\end{algorithmic} 
		\end{adjustwidth}
		\label{SLP_Radial}
	\end{algorithm}
	In general, constraints \eqref{opf_socr_RSOC_Slei} may not all be active at the optimum; thereby making the relaxation inexact and therefore infeasible in practice. In this case, Algorithm~\ref{SLP_SOCP} is extended to Algorithm~\ref{SLP_Radial}, which now includes supporting hyperplanes \emph{in conjunction with} supporting halfspaces of a set of the form \eqref{QoL} to recover an $\epsilon-$effective solution to Problem~\ref{opf_radial}. However, instead of directly using supporting hyperplanes of the form $\{ (w,w_{i}) | w_{i} = f_{ij}^{\rm a}(w,w^{(k)}) \} $ (to the nonconvex set of the form \eqref{QoL}), which can lead to infeasible LP problems, Algorithm~\ref{SLP_Radial} introduces non-negative slack variables $r_{ij}$ whose purpose is twofold. The first is to prevent infeasible LP problems, especially during the first few iterations of the algorithm when it is initialized from a poor-quality starting point. The second is to guide convergence by adaptively tuning the weight parameter $\rho_{ij}^{(k-1)}$. More specifically, when the condition on line 8 is true, $\rho_{ij}^{(k-1)}$ is increased by a factor of $\gamma$, as long as it remains smaller than a predefined upper limit $\overline{\rho}$. Set $\mathcal{I}_{ij}$ registers all previous iteration numbers that qualify as supporting halfspaces for \eqref{QoL}. The condition on line 5 is satisfied when constraint \eqref{QoL} is violated at iteration $k$, in which case the supporting halfspace $\{ (w,w_{i}) | w_{i} - f_{ij}^{\rm a}(w,w^{(k-1)}) \ge 0 \} $ is added to Problem \ref{SLP_radial} at $k+1$.
	
	In contrast to Algorithm~\ref{SLP_SOCP}, there is no theoretical guarantee that Algorithm~\ref{SLP_Radial} will converge to a feasible solution of Problem~\ref{opf_radial}, let alone an optimal one. In fact, \cite{Fletcher1998_SLPwithFilter} proves that an SLP algorithm in conjunction with a trust-region and an NLP filter applied to a general NLP with \emph{only} inequality constraints is guaranteed to converge to a feasible solution (global convergence\footnote{Global convergence of an algorithm entails convergence to a local optimum from \emph{any} starting point. Global convergence on general NLP problems should not be confused with convergence to a global optimum.}), from an arbitrary starting point that lies in the non-empty bounded region of the subset of the linear inequality constraints. However, the proof is predicated on the assumption that a Newton-like (second-order) feasibility restoration phase is used if an LP subproblem becomes infeasible, which means that LP solvers cannot be used exclusively. A geometric interpretation of Algorithm \ref{SLP_Radial} on a small nonconvex two-dimensional problem can be found in Appendix~\ref{slp_example}. The next section extends Algorithm~\ref{SLP_Radial} to meshed networks. 
	
\subsection{Proposed SLP algorithm for the OPF in meshed networks}\label{subsec_SLPforOPFMeshed}
	
	Problem \ref{opf_radial}, and therefore Algorithm~\ref{SLP_Radial}, is generally inexact on meshed networks as it lacks the treatment of the cycle constraints captured by \eqref{opf_alt_atan2}. The first-order Taylor series approximation of a function of the form $h_{ij}\left( w_{ij}^{{\rm i}},w_{ij}^{{\rm r}} \right) = {\rm atan2} \left( w_{ij}^{{\rm i}},w_{ij}^{{\rm r}} \right) $ at a point $\left(w_{ij}^{{\rm i},(k)},w_{ij}^{{\rm r},(k)} \right) $ can be written as
	\begin{multline}\label{atan2_Taylor}
		h_{ij}^{\rm a}\left(w_{ij}^{{\rm i}},w_{ij}^{{\rm r}},w_{ij}^{{\rm i},(k)},w_{ij}^{{\rm r},(k)} \right) = {\rm atan2} \left(w_{ij}^{{\rm i},(k)},w_{ij}^{{\rm r},(k)} \right) + \\
		\begin{bmatrix}
			\displaystyle \frac{-w_{ij}^{{\rm r},(k)}}{\left(w_{ij}^{{\rm r},(k)}\right)^2+\left(w_{ij}^{{\rm i},(k)}\right)^2} \\
			\displaystyle \frac{w_{ij}^{{\rm i},(k)}}{\left(w_{ij}^{{\rm r},(k)}\right)^2+\left(w_{ij}^{{\rm i},(k)}\right)^2}
		\end{bmatrix}^{T}
		\begin{bmatrix}
			w_{ij}^{{\rm i}} - w_{ij}^{{\rm i},(k)} \\
			w_{ij}^{{\rm r}} - w_{ij}^{{\rm r},(k)}
		\end{bmatrix}.
	\end{multline}
	By defining 
	\begin{align}
		H = \left[ \left( \theta_{i}^{(k)} - \theta_{j}^{(k)} - {\rm atan2} \left(w_{ij}^{{\rm i},(k)},w_{ij}^{{\rm r},(k)} \right) \right)_{ij \in \mathcal{L}} \right],
	\end{align}
	the SLP algorithm is described in Algorithm~\ref{SLP_Meshed}, which uses the same slack variables $r_{ij}$ in the Taylor series approximation of \eqref{opf_alt_atan2} to circumvent infeasible LP problems that may arise if \eqref{atan2_Taylor} is used directly instead of \eqref{opf_slp_atan1} and \eqref{opf_slp_atan2}.
	
	\begin{algorithm}[!t]
		\caption{SLP for the OPF in meshed networks.}
		\begin{adjustwidth}{-1em}{}
			\begin{algorithmic}[1]
				\scriptsize
				\STATE \parbox[t]{\dimexpr\linewidth-0.75cm}{\textbf{Initialization:}
					Same as in Algorithm~\ref{SLP_Radial}. \strut}
				\algrule
				\WHILE {$\max \left( F \cup H \right) > \epsilon$ \textbf{and} ${\rm max}\left( G \right) > \epsilon^{\rm th}$ \textbf{and} $ k \le \overline{k}$}
				\STATE \parbox[t]{\dimexpr\linewidth-0.75cm}{Obtain $\left( w_{i}^{(k)},w_{ij}^{{\rm r},(k)},w_{ij}^{{\rm i},(k)}, p_{ij}^{(k)}, q_{ij}^{(k)} \right) $ by solving
					\begin{subequations}\label{opf_slp}	
						\begin{align}
						& \hspace{-2cm} \underset {\substack{p_{gi},p_{g},\alpha_{n},q_{gi},w_{i},\theta_{i}\\w_{ij}^{\rm r},w_{ij}^{\rm i},p_{ij},q_{ij},r_{ij}}} 
						{\mbox{ minimize}} \eqref{opf_socr_slp_objective} + \sum_{ij \in \mathcal{L}} \rho_{ij}^{(k-1)}r_{ij} \hspace{-3cm} & & \label{opf_slp_objective} \\
						& \hspace{-1.5cm} \text{ subject to \cref{SLP_radial_shared,SLP_radial_cut0,SLP_radial_cut1,SLP_radial_cut_thermal,opf_slp_r}} \hspace{-2cm} & & \label{opf_slp_shared} \\
						& \hspace{-1.25cm}h_{ij}^{\rm a}\left(w_{ij}^{{\rm i}},w_{ij}^{{\rm r}},w_{ij}^{{\rm i},(k-1)},w_{ij}^{{\rm r},(k-1)} \right) \le r_{ij}, \hspace{-1.5cm} & ij \in \mathcal{L}, & \label{opf_slp_atan1} \\
						& \hspace{-1.25cm} h_{ij}^{\rm a}\left(w_{ij}^{{\rm i}},w_{ij}^{{\rm r}},w_{ij}^{{\rm i},(k-1)},w_{ij}^{{\rm r},(k-1)} \right) \ge -r_{ij}, \hspace{-1.5cm} & ij \in \mathcal{L}, & \label{opf_slp_atan2}
						\end{align}
					\end{subequations} \strut} 
				\STATE Lines 4 to 16 in Algorithm~\ref{SLP_Radial}.
				\STATE {$k:= k + 1 $.}
				\ENDWHILE
			\end{algorithmic} 
		\end{adjustwidth}
		\label{SLP_Meshed}
	\end{algorithm} 

	The weight parameter is \emph{automatically} tuned exactly as described in Algorithm~\ref{SLP_Radial}, and numerical experiments have shown that setting $\rho_{ij}^{(0)}$ to a fixed high value, say six orders of magnitude larger than $ \max\left( c_{2,gi},c_{1,gi} \right) $ would result in fast convergence to a feasible but poor-quality solution. Instead, setting $\rho_{ij}^{(0)}$ to around one order of magnitude larger than $ \max\left( c_{2,gi},c_{1,gi} \right) $ and adaptively increasing it as described on lines 8 to 10 in Algorithm~\ref{SLP_Radial} was numerically found to strike a good tradeoff between computation time and solution quality.
	
	Aside from the tolerances and $\zeta$, the convergence of Algorithm~\ref{SLP_Meshed} hinges on a \emph{single} parameter $\rho_{ij}$, which is then automatically tuned by introducing parameters $\gamma$ and $\overline{\rho}$. More importantly, an AC feasibility restoration phase is not needed in Algorithm~\ref{SLP_Meshed} owing to the iteratively refined polyhedral outer approximation of (i) the nonconvex set of the form \eqref{QoL} by the supporting halfspaces in \eqref{SLP_radial_cut1} in conjunction with the supporting hyperplanes in \eqref{SLP_radial_cut0}, and (ii) the convex set of the form $\{ (p_{ij},q_{i}) | p_{ij}^{2} + q_{ij}^{2} \le \overline{s}^{2}_{ij}, ij \in \mathcal{L} \cup \mathcal{L}_{\rm t} \} $ by the supporting halfspaces in \eqref{SLP_radial_cut_thermal}. Extensive numerical evaluation of Algorithms~\ref{SLP_Radial} and~\ref{SLP_Meshed} is conducted in the next section to demonstrates their exactness, computational efficiency, and robustness against the choice of starting point.


\section{Numerical evaluation}\label{sec_numericalevalulation} 
	
	In this experimental setup, Julia v1.4.0 \cite{Bezanson2017_Julia} is used as a programming language along with JuMP v0.21.1 \cite{DunningHuchetteLubin2017_JuMP} as a frontend modeling language for the optimization problems. The original NLP problem in \eqref{opf_alt} is solved using KNITRO v12.2.2 \cite{KNITRO} with default settings. The LP problems in Algorithm~\ref{SLP_SOCP} and Algorithm~\ref{SLP_Meshed} are solved using Gurobi v9.0.2 \cite{Gurobi2019} with the Barrier method (``Method = 2'') and the following parameters:
	\begin{itemize}
		\item ``Crossover = 0'', which disables crossover,\footnote{The crossover strategy transforms the interior solution returned by the Barrier method into a basic solution to be used at the root node of an MILP, and is therefore not needed for the continuous models in this work.}
		\item ``ObjScale = -1'', which uses the reciprocal of the maximum coefficient as scaling for the objective function,
		\item ``ScaleFlag = -1'', which implements model scaling for improving the numerical properties of the constraint matrix, 
		\item ``Presolve = 2'', which implements an aggressive presolve strategy in an effort to tighten the model.
	\end{itemize}
	For all the test cases, Algorithms~\ref{SLP_Radial} and~\ref{SLP_Meshed} are initialized with $\epsilon = 10^{-5}$, $\epsilon^{\rm th} = 10^{-3}$, $\zeta = 0.9$, $\rho_{ij}^{(0)} = 10 \max \left( \left( {\max}\left( c_{2,gi},c_{1,gi}\right)  \right) _{gi \in \mathcal{G}} \right)$, $\gamma = 5$, and $\overline{\rho}=\gamma^{4}\rho_{ij}^{(0)}$. Finally, let the \emph{relative} optimality gap be defined as
	\begin{align}\label{eq_optimality}
		Gap=\left( \left( P^{\dagger}_{\rm IPM}- P^{\dagger}_{ \rm SLP}\right) /P^{\dagger}_{ \rm IPM}\right) \times 100,
	\end{align}
	where $P^{\dagger}_{\rm IPM}$ is the locally optimal solution obtained by KNITRO and $P^{\dagger}_{\rm SLP}$ is the $\epsilon$-effective solution obtained by Algorithm~\ref{SLP_Radial} or Algorithm~\ref{SLP_Meshed}.

\subsection{Exactness}\label{sec_optandfeasbility}
	\subsubsection{Radial networks}
	The radial test cases consist of the 33-node \cite{Baran1989_33bus}, 69-node \cite{Savier2007_69bus}, and 119-node \cite{Zhang2007_119bus} IEEE distribution systems with the tie lines removed to ensure a radial topology.\footnote{It is worth noting that keeping the tie lines results in weakly meshed topologies which necessitate the use of Algorithm~\ref{SLP_Meshed} instead of Algorithm~\ref{SLP_Radial}.} The exactness and computational performance of Algorithm~\ref{SLP_Radial} for three different starting points is shown in Table~\ref{Table_OPF_Radial}. More specifically, let FS1 designate a flat start consisting of $(v_{i}^{(0)}=1,\theta_{i}^{(0)}=0 )$, FS2 designate a flat start consisting of $(v_{i}^{(0)}=\underline{v}_{i},\theta_{i}^{(0)}=0 )$, and FS3 designate a flat start consisting of $(v_{i}^{(0)}=\overline{v}_{i},\theta_{i}^{(0)}=0 )$, for all $i \in \mathcal{B}$.
	\begin{table}[!t]
		\centering
		\renewcommand{\arraystretch}{1.3}
		\caption{Exactness and computational performance of Algorithm~\ref{SLP_Radial} compared to KNITRO on 33-node, 69-node, and 119-node IEEE distribution systems for three different starting points.}
		\resizebox{\linewidth}{!}{%
			\begin{tabular}{| l | c | c | c | c | c | c |} 	\hline
				\multicolumn{1}{|l|}{\textbf{Test}} & \multicolumn{2}{c|}{\textbf{NLP (KNITRO)}} & \multicolumn{4}{c|}{\textbf{Algorithm~\ref{SLP_Radial} (Gurobi)}} 
				\\\cline{2-7}
				\textbf{Case} & \textbf{Cost} ($\SI{}{\$}$) & \textbf{Time} ($\SI{}{\second}$) & \textbf{Time} ($\SI{}{\second}$) & $\boldsymbol{k}$ & $\boldsymbol{{\rm mean}(F)}$ & $\boldsymbol{Gap}$ (\%) \\\hline
				\multicolumn{7}{|c|}{\textbf{FS1}} \\\hline													
				33bus	&	386.03	&	0.05	&	0.01	&	2	&	2.60E-07	&	6.65E-04	\\\hline
				69bus	&	434.18	&	0.08	&	0.03	&	3	&	1.35E-10	&	-2.80E-04	\\\hline
				119bus	&	2434.65	&	0.08	&	0.04	&	3	&	7.96E-11	&	8.17E-06	\\\hline
				\multicolumn{7}{|c|}{\textbf{FS2}} \\\hline													
				33bus	&	386.03	&	0.05	&	0.01	&	2	&	2.53E-07	&	6.45E-04	\\\hline
				69bus	&	434.18	&	0.07	&	0.04	&	3	&	1.58E-10	&	9.93E-06	\\\hline
				119bus	&	2434.65	&	0.08	&	0.04	&	3	&	1.12E-10	&	-1.01E-05	\\\hline
				\multicolumn{7}{|c|}{\textbf{FS3}} \\\hline													
				33bus	&	386.03	&	0.06	&	0.01	&	2	&	2.50E-09	&	6.58E-04	\\\hline
				69bus	&	434.18	&	0.07	&	0.03	&	3	&	1.37E-10	&	-1.80E-04	\\\hline
				119bus	&	2434.65	&	0.07	&	0.04	&	3	&	1.12E-12	&	-1.04E-05	\\\hline
		\end{tabular}}
		\label{Table_OPF_Radial}
	\end{table}
	In Table~\ref{Table_OPF_Radial}, columns 2 and 3 show the objective function value \eqref{opf_alt_objective} at the optimum and the associated computation time of KNITRO, respectively. The computation times and the associated number of iterations of Algorithm~\ref{SLP_Radial} are listed in columns 4 and 5, respectively. Column 6 shows the mean violations of constraints \eqref{opf_alt_RSOC}, and the last column shows the relative optimality gap of Algorithm~\ref{SLP_Radial} defined in \eqref{eq_optimality}. It is evident from Table~\ref{Table_OPF_Radial} that Algorithm~\ref{SLP_Radial} converges on all the three test cases from all three starting point strategies to high-quality solutions in computation times within the same order of magnitude as KNITRO's. In particular, the mean violations of constraints \eqref{opf_alt_RSOC} and the (relative) optimality gaps are on average around $10^{-8}$ and $10^{-4}$\%, respectively. In fact, Algorithm~\ref{SLP_Radial} converges from \emph{any} starting point such that $\underline{v}_{i} \le v_{i}^{(0)} \le \overline{v}_{i}$ for all $i \in \mathcal{B}$. The robustness of Algorithm~\ref{SLP_Radial} against the choice of starting point will be analyzed in more detail in Section~\ref{sec_Robustness}.\footnote{Please note that Algorithm 3 is for general networks and can thus be applied to any grid topology, including radial and weakly meshed topologies. However, since radial networks are acyclic, the hyperplanes \eqref{atan2_Taylor} approximating the cycle constraints are redundant and keeping them may add unnecessary computational burden. It is for this reason that Algorithm~\ref{SLP_Radial} is used instead of Algorithm~\ref{SLP_Meshed} for radial networks.}
	
	\subsubsection{Meshed networks}
	Algorithm \ref{SLP_Meshed} is implemented on a wide range of meshed network test cases, namely, 27 instances available with MATPOWER \cite{MATPOWER}, and 108 of the more difficult instances from PGLib-OPF \cite{Babaeinejadsarookolaee2019_PGlibOPF}, for a total of 135 test cases with up the 3375 buses. The PGLib-OPF test cases are designed specifically to incorporate key network parameters such as line thermal limit and small angle difference, which are critical in assessing the robustness of OPF optimization algorithms. In particular, active power increase (API) cases are designed to emulate heavily loaded systems (with binding thermal limit constraints) whereas the small angle difference (SAD) cases are designed to emulate small angle differences in some practical systems.\footnote{The angle difference in the PGLib-OPF SAD are $\pm 10^{\circ}$.} The exactness and computational performance of Algorithm~\ref{SLP_Meshed} are demonstrated on four different choices of starting points, namely three flat starts and a DC OPF warm start. The three flat starts are the same ones defined in the previous section, i.e., FS1, FS2, and FS3. The exactness and computational performance of Algorithm~\ref{SLP_Meshed} compared to KNITRO's on MATPOWER and PGLib-OPF test cases with sizes ranging from 2383 to 3375 buses for FS1 and a DC OPF warm start are shown in Table~\ref{Table_OPF_FS1} and Table~\ref{Table_OPF_DC}, respectively. The DC OPF is infeasible on most PGLib-OPF SAD test cases (also evidenced in \cite{Babaeinejadsarookolaee2019_PGlibOPF}), which is why Table~\ref{Table_OPF_DC} only shows the results for the corresponding MATPOWER, PGLib-OPF TYP, and PGLib-OPF API test cases. For lack of space, not all results for the 135 test cases are shown here. Results for the remaining test cases under FS1 and the DC OPF warm start as well as those under FS2 and FS3 can be found in \cite{SLP_Results}. Numerical evaluation of Algorithm~\ref{SLP_SOCP} under FS1 can also be found in \cite{SLP_Results}. In Tables~\ref{Table_OPF_FS1} and~\ref{Table_OPF_DC}, columns 2 and 3 show the objective function value \eqref{opf_alt_objective} at the optimum and the associated computation time of KNITRO, respectively. The computation times and the associated number of iterations of Algorithm~\ref{SLP_Meshed} are listed in columns 4 and 5, respectively. Column 6 shows the mean violations of constraints \eqref{opf_alt_RSOC} and \eqref{opf_alt_atan2}, respectively, and the last column shows the \emph{relative} optimality gap defined in \eqref{eq_optimality}.
	
	As evidenced in \cite{SLP_Results} and Tables~\ref{Table_OPF_FS1} and~\ref{Table_OPF_DC}, Algorithm~\ref{SLP_Meshed} converges on all the 135 test cases from all four starting point strategies to high-quality solutions in computation times generally within the same order of magnitude as KNITRO's. In particular, the mean violations of constraints \eqref{opf_alt_RSOC} and \eqref{opf_alt_atan2} and the (relative) optimality gaps are on average around $10^{-7}$ and $10^{-3}$\%, respectively. In fact, Algorithm~\ref{SLP_Meshed} converges for \emph{any} starting point such that $\underline{v}_{i} \le v_{i}^{(0)} \le \overline{v}_{i}$ for all $i \in \mathcal{B}$ and $\theta_{i}^{(0)}$ such that $\underline{\theta}_{ij} \le \theta_{ij}^{(0)} \le \overline{\theta}_{ij}$ for all $ij \in \mathcal{L}$, both of which are trivial to ensure in practice. The robustness of Algorithm~\ref{SLP_Meshed} against the choice of starting point will be further analyzed in the next section. More interestingly, most of the cases that take more than 10 iterations to converge under the generic parameter setting described above can converge faster under a different value of $\rho_{ij}^{(0)}$. For instance, MATPOWER's Polish test case 2383wp with $\rho_{ij}^{(0)} = 100\max\left( c_{2,gi},c_{1,gi} \right)$ converges in 6 iterations and 9 seconds and 8 iterations and 12 seconds with the DC Warm start and FS1, respectively. A similar performance enhancement is witnessed on the PGLib-OPF TYP and PGLib-OPF SAD versions of the 2383wp test case. Another example is PGLib-OPF SAD's case2316\_sdet which converges in 10 iterations and 18 seconds under FS1. In practice, a power system operator can tailor $\rho_{ij}^{(0)}$ and $\gamma$ specifically to the power system at hand for an ideal performance. Furthermore, as mentioned earlier, not all the supporting half spaces in \eqref{SLP_radial_cut1} are necessary for convergence, which means that previous inactive ones can be removed in favor of faster convergence.
	\begin{table}[!t]
		\centering
		\renewcommand{\arraystretch}{1.3}
		\caption{Exactness and computational performance of Algorithm~\ref{SLP_Meshed} compared to KNITRO on all the MATPOWER and PGLib-OPF test cases with sizes ranging from 2383 to 3375 buses for starting point FS1.}
		\resizebox{\linewidth}{!}{%
			\begin{tabular}{| l | c | c | c | c | c | c |} 	\hline
				\multicolumn{1}{|l|}{\textbf{Test}} & \multicolumn{2}{c|}{\textbf{NLP (KNITRO)}} & \multicolumn{4}{c|}{\textbf{Algorithm~\ref{SLP_Meshed} (Gurobi)}} 
				\\\cline{2-7}
				\textbf{Case} & \textbf{Cost} ($\SI{}{\$}$) & \textbf{Time} ($\SI{}{\second}$) & \textbf{Time} ($\SI{}{\second}$) & $\boldsymbol{k}$ & $\boldsymbol{{\rm mean}(F \cup H)}$ & $\boldsymbol{Gap}$ (\%) \\\hline
				\multicolumn{7}{|c|}{\textbf{MATPOWER}} \\\hline
				2383	&	1868511.83	&	3.8	&	37.4	&	15	&	1.1E-07	&	-9.1E-03	\\\hline
				2736	&	1307883.13	&	3.5	&	5.5	&	4	&	8.5E-07	&	1.7E-03	\\\hline
				2737	&	777629.30	&	3.3	&	4.8	&	4	&	1.2E-07	&	-4.1E-05	\\\hline
				2746wop	&	1208279.81	&	2.9	&	5.1	&	4	&	1.1E-06	&	3.2E-03	\\\hline
				2746wp	&	1631775.10	&	2.9	&	4.9	&	4	&	7.5E-07	&	2.9E-03	\\\hline
				2848	&	53022.04	&	19.6	&	26.8	&	12	&	4.3E-07	&	5.8E-05	\\\hline
				2868	&	79794.68	&	7.0	&	22.7	&	11	&	1.6E-06	&	-2.6E-04	\\\hline
				2869	&	133999.32	&	7.0	&	25.7	&	9	&	8.9E-07	&	3.4E-04	\\\hline
				3012	&	2591706.57	&	4.1	&	9.3	&	6	&	8.4E-08	&	1.2E-04	\\\hline
				3120	&	2142703.77	&	4.0	&	9.4	&	6	&	6.0E-08	&	2.3E-04	\\\hline
				3375	&	7412030.68	&	6.9	&	10.8	&	6	&	1.1E-06	&	4.7E-05	\\\hline
				\multicolumn{7}{|c|}{\textbf{PGLib-OPF TYP}} \\\hline													
				2000	&	1228487.06	&	2.7	&	7.1	&	5	&	2.0E-06	&	5.4E-04	\\\hline
				2316	&	1775325.56	&	2.9	&	33.6	&	15	&	3.1E-07	&	-2.4E-03	\\\hline
				2383	&	1868191.64	&	3.9	&	36.2	&	15	&	1.4E-07	&	-1.1E-02	\\\hline
				2736	&	1308015.00	&	3.3	&	5.1	&	4	&	7.5E-07	&	2.0E-03	\\\hline
				2737	&	777727.68	&	2.7	&	5.0	&	4	&	1.6E-07	&	4.8E-04	\\\hline
				2746wop	&	1208258.50	&	2.8	&	5.0	&	4	&	1.2E-06	&	2.7E-03	\\\hline
				2746wp	&	1631707.93	&	3.0	&	5.3	&	4	&	7.4E-07	&	2.8E-03	\\\hline
				2848	&	1286608.21	&	6.8	&	18.4	&	10	&	3.8E-07	&	-6.1E-04	\\\hline
				2853	&	2052386.73	&	4.9	&	29.4	&	11	&	1.0E-07	&	-4.4E-04	\\\hline
				2868	&	2009605.33	&	6.8	&	17.9	&	10	&	6.9E-07	&	-1.8E-04	\\\hline
				2869	&	2462790.45	&	6.7	&	50.1	&	13	&	5.4E-08	&	-5.7E-04	\\\hline
				3012	&	2600842.77	&	4.6	&	9.6	&	6	&	4.9E-08	&	3.3E-04	\\\hline
				3120	&	2147969.11	&	4.2	&	15.5	&	8	&	2.1E-08	&	-4.5E-04	\\\hline
				3375	&	7438169.48	&	6.2	&	21.9	&	10	&	1.7E-07	&	-4.6E-03	\\\hline
				\multicolumn{7}{|c|}{\textbf{PGLib-OPF API}} \\\hline													
				2000	&	1285030.37	&	15.1	&	35.9	&	11	&	1.2E-07	&	-1.4E-04	\\\hline
				2316	&	2189027.77	&	3.8	&	41.4	&	16	&	9.4E-08	&	-2.0E-03	\\\hline
				2383	&	279125.83	&	1.8	&	6.7	&	6	&	8.2E-06	&	6.2E-10	\\\hline
				2736	&	653940.01	&	3.6	&	11.0	&	7	&	3.2E-07	&	-4.5E-03	\\\hline
				2737	&	369197.40	&	4.1	&	13.9	&	8	&	1.8E-07	&	8.5E-04	\\\hline
				2746wop	&	511658.63	&	1.8	&	7.0	&	7	&	1.8E-06	&	-7.1E-10	\\\hline
				2746wp	&	581825.06	&	1.9	&	5.4	&	6	&	5.4E-06	&	3.3E-09	\\\hline
				2848	&	1475997.12	&	5.9	&	46.1	&	16	&	4.0E-07	&	-2.3E-03	\\\hline
				2853	&	2457850.66	&	6.2	&	39.6	&	13	&	4.6E-07	&	-7.2E-04	\\\hline
				2868	&	2340496.19	&	348.0	&	39.0	&	15	&	2.5E-07	&	-3.8E-04	\\\hline
				2869	&	2934160.74	&	6.4	&	46.8	&	13	&	9.9E-08	&	-1.3E-02	\\\hline
				3012	&	728874.00	&	2.7	&	18.4	&	6	&	1.1E-06	&	-1.0E-04	\\\hline
				3120	&	984280.13	&	6.3	&	39.1	&	15	&	4.1E-08	&	-5.3E-03	\\\hline
				3375	&	5847780.45	&	6.6	&	17.4	&	8	&	8.3E-07	&	-3.2E-03	\\\hline
				\multicolumn{7}{|c|}{\textbf{PGLib-OPF SAD}} \\\hline													
				2000	&	1230320.61	&	3.8	&	31.5	&	12	&	1.9E-07	&	-6.0E-05	\\\hline
				2316	&	1775330.54	&	3.5	&	40.2	&	16	&	4.5E-07	&	-3.9E-03	\\\hline
				2383	&	1912693.95	&	4.4	&	63.7	&	21	&	1.6E-07	&	-3.7E-02	\\\hline
				2736	&	1327303.04	&	4.6	&	4.1	&	3	&	3.0E-06	&	1.8E-02	\\\hline
				2737	&	791525.57	&	4.5	&	6.9	&	4	&	2.1E-07	&	3.4E-04	\\\hline
				2746wop	&	1234337.97	&	3.4	&	4.1	&	3	&	2.8E-06	&	1.2E-02	\\\hline
				2746wp	&	1667583.13	&	4.7	&	5.5	&	4	&	1.4E-06	&	4.1E-03	\\\hline
				2848	&	1289004.65	&	12.7	&	44.8	&	18	&	1.1E-07	&	-3.9E-03	\\\hline
				2853	&	2070064.41	&	5.9	&	44.2	&	14	&	2.3E-07	&	-1.2E-03	\\\hline
				2868	&	2022411.38	&	7.1	&	23.5	&	12	&	4.6E-07	&	-5.7E-04	\\\hline
				2869	&	2468863.64	&	8.2	&	48.1	&	12	&	4.3E-08	&	-1.3E-04	\\\hline
				3012	&	2621315.02	&	5.0	&	11.7	&	6	&	2.5E-08	&	4.2E-04	\\\hline
				3120	&	2178216.10	&	6.7	&	50.8	&	17	&	6.3E-08	&	-2.3E-03	\\\hline
				3375	&	7438169.48	&	6.8	&	21.9	&	10	&	6.7E-08	&	-1.1E-03	\\\hline
		\end{tabular}}
		\label{Table_OPF_FS1}
	\end{table}
	\begin{table}[!t]
		\centering
		\renewcommand{\arraystretch}{1.3}
		\caption[Performance]{Exactness and computational performance of Algorithm~\ref{SLP_Meshed} compared to KNITRO on MATPOWER, PGLib-OPF TYP\footnotemark, and PGLib-OPF API test cases with sizes ranging from 2383 to 3375 buses for a DC OPF warm start.}
		\resizebox{\linewidth}{!}{%
			\begin{tabular}{| l | c | c | c | c | c | c |} 	\hline
				\multicolumn{1}{|l|}{\textbf{Test}} & \multicolumn{2}{c|}{\textbf{NLP (KNITRO)}} & \multicolumn{4}{c|}{\textbf{Algorithm~\ref{SLP_Meshed} (Gurobi)}} 
				\\\cline{2-7}
				\textbf{Case} & \textbf{Cost} ($\SI{}{\$}$) & \textbf{Time} ($\SI{}{\second}$) & \textbf{Time} ($\SI{}{\second}$) & $\boldsymbol{k}$ & $\boldsymbol{{\rm mean}(F \cup H)}$ & $\boldsymbol{Gap}$ (\%) \\\hline
				\multicolumn{7}{|c|}{\textbf{MATPOWER}} \\\hline
				2383	&	1868511.83	&	3.6	&	51.3	&	18	&	1.1E-07	&	2.8E-02	\\\hline
				2736	&	1307883.13	&	3.1	&	3.6	&	3	&	1.3E-06	&	5.6E-03	\\\hline
				2737	&	777629.30	&	4.7	&	5.1	&	4	&	3.5E-07	&	4.3E-04	\\\hline
				2746wop	&	1208279.81	&	3.3	&	3.6	&	3	&	3.0E-06	&	1.6E-02	\\\hline
				2746wp	&	1631775.10	&	3.5	&	3.6	&	3	&	4.9E-06	&	2.0E-02	\\\hline
				2848	&	53022.12	&	26.6	&	18.0	&	9	&	2.5E-06	&	-2.2E-04	\\\hline
				2868	&	79794.68	&	9.7	&	27.9	&	12	&	1.2E-06	&	-1.9E-04	\\\hline
				2869	&	133999.31	&	5.1	&	36.0	&	11	&	1.3E-07	&	1.0E-05	\\\hline
				3012	&	2591706.57	&	5.1	&	9.0	&	6	&	2.1E-07	&	1.1E-03	\\\hline
				3120	&	2142703.77	&	5.2	&	5.5	&	4	&	3.5E-07	&	2.9E-04	\\\hline
				3375	&	7412030.68	&	6.3	&	12.7	&	7	&	5.5E-07	&	-6.3E-03	\\\hline
				\multicolumn{7}{|c|}{\textbf{PGLib-OPF TYP}} \\\hline													
				2000	&	1228487.06	&	2.8	&	8.9	&	6	&	6.4E-07	&	8.4E-05	\\\hline
				2316	&	1775325.56	&	2.8	&	21.7	&	11	&	1.1E-06	&	-4.5E-03	\\\hline
				2383	&	1868191.64	&	3.7	&	62.4	&	22	&	3.1E-06	&	1.8E-02	\\\hline
				2736	&	1308015.00	&	2.8	&	3.6	&	3	&	1.2E-06	&	5.9E-03	\\\hline
				2737	&	777727.68	&	2.5	&	5.2	&	4	&	3.6E-07	&	4.1E-04	\\\hline
				2746wop	&	1208258.50	&	2.6	&	3.7	&	3	&	3.1E-06	&	1.8E-02	\\\hline
				2746wp	&	1631707.93	&	2.8	&	3.6	&	3	&	4.8E-06	&	1.9E-02	\\\hline
				2848	&	1286608.21	&	12.3	&	14.6	&	8	&	3.1E-07	&	-5.9E-04	\\\hline
				2853	&	2052386.73	&	4.9	&	34.6	&	12	&	1.3E-07	&	-2.1E-03	\\\hline
				2868	&	2009605.33	&	6.4	&	14.2	&	9	&	1.9E-07	&	-1.3E-04	\\\hline
				2869	&	2462790.45	&	5.6	&	38.0	&	11	&	1.8E-07	&	-2.1E-04	\\\hline
				3012	&	2600842.77	&	4.2	&	13.6	&	7	&	4.7E-08	&	2.2E-04	\\\hline
				3120	&	2147969.11	&	7.7	&	19.4	&	10	&	1.4E-08	&	-7.8E-05	\\\hline
				3375	&	7438169.48	&	6.5	&	26.6	&	11	&	7.3E-08	&	-1.8E-03	\\\hline
				\multicolumn{7}{|c|}{\textbf{PGLib-OPF API}} \\\hline													
				2000	&	1285030.37	&	17.3	&	28.3	&	10	&	1.8E-07	&	-4.5E-04	\\\hline
				2316	&	2189027.77	&	3.3	&	33.2	&	9	&	1.1E-06	&	-2.6E-03	\\\hline
				2383	&	279125.83	&	2.3	&	7.3	&	6	&	8.7E-06	&	6.2E-10	\\\hline
				2736	&	653940.01	&	3.7	&	10.4	&	7	&	3.3E-07	&	-5.1E-03	\\\hline
				2737	&	369197.40	&	3.8	&	16.7	&	9	&	2.6E-08	&	-1.7E-03	\\\hline
				2746wop	&	511658.63	&	1.3	&	17.1	&	10	&	2.6E-06	&	-9.7E-06	\\\hline
				2746wp	&	581825.06	&	1.8	&	10.0	&	8	&	2.6E-06	&	1.1E-08	\\\hline
				2848	&	1475997.12	&	5.2	&	49.9	&	18	&	4.7E-07	&	-1.7E-03	\\\hline
				2853	&	2457850.66	&	6.9	&	58.3	&	16	&	9.6E-08	&	-4.0E-04	\\\hline
				2868	&	2340496.19	&	5.4	&	46.9	&	18	&	4.9E-07	&	-4.2E-04	\\\hline
				2869	&	2934160.74	&	7.1	&	44.2	&	12	&	7.0E-08	&	-8.0E-03	\\\hline
				3012	&	728874.00	&	2.7	&	14.0	&	5	&	1.6E-06	&	-2.0E-05	\\\hline
				3120	&	984280.13	&	5.5	&	37.9	&	15	&	5.3E-08	&	-1.3E-02	\\\hline
				3375	&	5847780.45	&	5.9	&	26.9	&	10	&	2.4E-07	&	-5.6E-04	\\\hline
		\end{tabular}}
		\label{Table_OPF_DC}
	\end{table}
	\footnotetext{TYP is short for Typical Operating Conditions in \cite{Babaeinejadsarookolaee2019_PGlibOPF}.}
	
\subsection{Robustness}\label{sec_Robustness}

\begin{figure}[t!]
		\centering{
			\includegraphics[width=1.0\columnwidth] {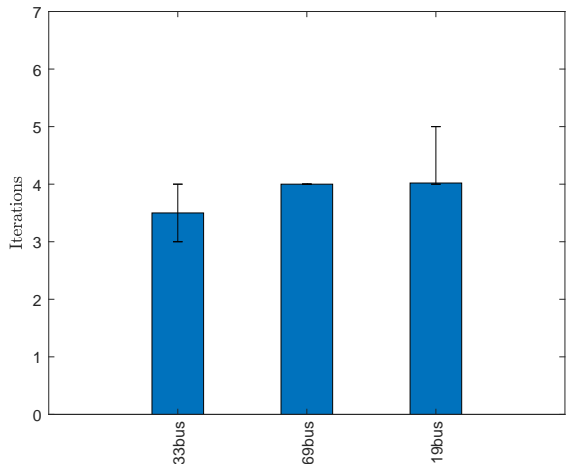}}
		\caption{Average number of iterations to convergence of Algorithm~\ref{SLP_Radial} on the 33-node, 69-node, and 119-node IEEE distribution systems. The error bars show the maximum and minimum number of iterations to convergence for each test case.}
		\label{fig_Radial_robustness}
	\end{figure}

\begin{figure}[t!]
		\centering{
			\includegraphics[width=1.0\columnwidth] {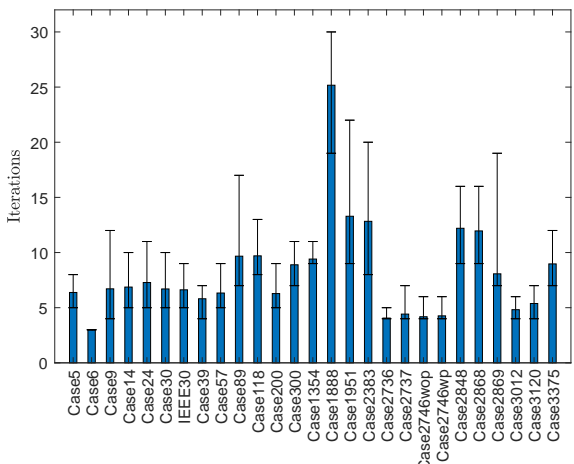}}
		\caption{Average number of iterations to convergence of Algorithm~\ref{SLP_Meshed} on the 27 meshed cases from MATPOWER \cite{MATPOWER}. The error bars show the maximum and minimum number of iterations to convergence for each test case.}
		\label{fig_MATPOWER_robustness}
	\end{figure}

	This section further analyses the robustness of Algorithms~\ref{SLP_Radial} and~\ref{SLP_Meshed} beyond the three flat starts identified in the previous section. Towards this aim, Algorithms~\ref{SLP_Radial} and~\ref{SLP_Meshed} are initialized from 100 different starting points with bus voltages $v_{i}^{(0)}$ drawn from a uniform distribution in the interval $\left( \underline{v}_{i}, \overline{v}_{i} \right) $, \emph{i.e.}, $v_{i}^{(0)} \sim \mathcal{U} \left( \underline{v}_{i}, \overline{v}_{i} \right) $ and $\theta_{i}^{(0)} = 0$ for all $i \in \mathcal{B}$. For each test case, Algorithms~\ref{SLP_Radial} and~\ref{SLP_Meshed} converge from \emph{all} the 100 different starting points to the same objective values obtained from KNITRO \cite{KNITRO}, and the average number of iterations to convergence is shown in Fig. \ref{fig_Radial_robustness} and Fig. \ref{fig_MATPOWER_robustness} for the three radial test cases and the 27 meshed cases from MATPOWER \cite{MATPOWER}, respectively. The error bars in Fig. \ref{fig_Radial_robustness} and Fig. \ref{fig_MATPOWER_robustness} show the maximum and minimum number of iterations to convergence for each test case. The robustness is corroborated by the following two observations: (i) for each test case, the two algorithms converge from \emph{all} the 100 different starting points to the same objective values obtained from KNITRO \cite{KNITRO}, and (ii) the average number of iterations to convergence for each test case is very close, if not equal, to the one obtained from starting point FS1 where all the bus voltages $v_{i}^{(0)}$ are set to 1 pu.
	A similar observation is witnessed on the remaining test cases from PGLib-OPF \cite{Babaeinejadsarookolaee2019_PGlibOPF}. Those remaining results are omitted here for ease of exposition but can be found in \cite{SLP_Results}. 

\subsection{LMP and Q-LMP}\label{sec_LMP}
	
\begin{figure}[t!]
		\centering
		\begin{subfigure}{\columnwidth}
			\psfrag{Bus}{\footnotesize Bus \normalsize}
			\psfrag{LMP}{\footnotesize LMP $\left( \SI{}{\$ \per \mega \watt} \right) $ \normalsize}
			\psfrag{AC OPF}{\footnotesize KNITRO \normalsize}
			\psfrag{Algorithm}{\footnotesize Algorithm~\ref{SLP_Meshed} \normalsize}
			\includegraphics[width=\columnwidth] {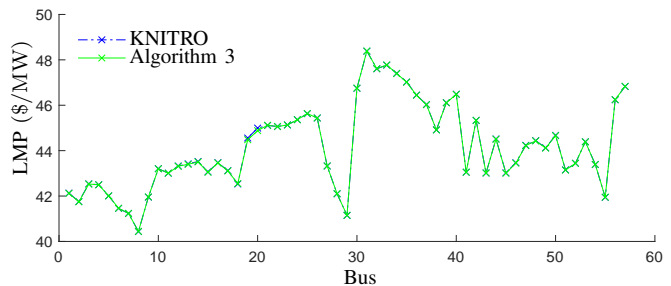}
			\caption{LMP.}
			\label{fig_LMP57}
		\end{subfigure}
		\begin{subfigure}{\columnwidth}
			\psfrag{Bus}{\footnotesize Bus \normalsize}
			\psfrag{QLMP}{\footnotesize Q-LMP $\left( \SI{}{\$ \per \mega \VAr} \right) $ \normalsize}
			\psfrag{AC OPF}{\footnotesize KNITRO \normalsize}
			\psfrag{Algorithm}{\footnotesize Algorithm~\ref{SLP_Meshed} \normalsize}
			\includegraphics[width=\columnwidth] {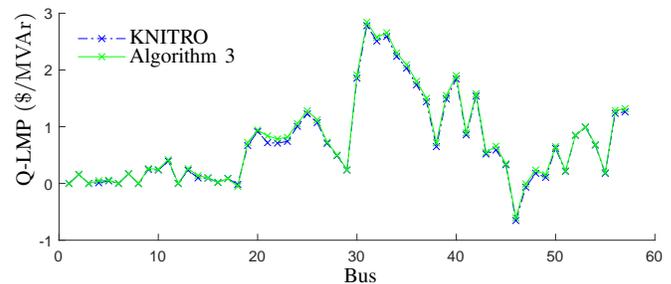}
			\caption{Q-LMP.}
			\label{fig_QLMP57}
		\end{subfigure}
		\caption{LMP and Q-LMP of KNITRO \cite{KNITRO} and Algorithm~\ref{SLP_Meshed} on the IEEE 57-bus test system \cite{MATPOWER}.}
		\label{fig_MP57}
	\end{figure}

	The exactness of the solutions from Algorithm~\ref{SLP_Meshed} also translates to accurate LMP and Q-LMP as shown in Table \ref{Table_LMP}. Table \ref{Table_LMP} also shows the maximum errors in nodal voltages and branch power flows. Table \ref{Table_LMP} affirms that Algorithm~\ref{SLP_Meshed} consistently produces accurate LMP and Q-LMP as well as accurate nodal voltages and branch power flows. This property entails that market operators already using LP solvers can not only generate more accurate LMP by using Algorithm~\ref{SLP_Meshed} compared to the classical DC OPF, but can now offer a new avenue for pricing reactive power through Q-LMP, without changing their LP solvers.
	\begin{table}[!t]
		\centering
		\renewcommand{\arraystretch}{1.3}
		\caption{Algorithm~\ref{SLP_Meshed}'s LMP, Q-LMP, node voltages, and branch power flows compared to those from KNITRO \cite{KNITRO}.}
		\resizebox{\linewidth}{!}{%
			\begin{tabular}{l c c c c} 	\hline \hline
				\multicolumn{1}{l}{\textbf{Test}} & \multicolumn{1}{c}{\textbf{Mean LMP}} & \multicolumn{1}{c}{\textbf{Mean Q-LMP}} & \multicolumn{1}{c}{\textbf{Max Voltage}} & \multicolumn{1}{c}{\textbf{Max $\boldsymbol{p_{ij}}$}}
				\\
				\textbf{Case} & \textbf{error} ($\SI{}{\$ \per \mega \VA}$) & \textbf{error} ($\SI{}{\$ \per \mega \VAr}$) & \textbf{error} (pu) & \textbf{error} (pu) \\\hline
				5	&	7.44E-05	&	8.80E-04	&	9.40E-06	&	7.86E-06	\\\hline
				6	&	8.06E-03	&	7.37E-03	&	4.10E-06	&	9.52E-04	\\\hline
				9	&	2.29E-03	&	1.84E-04	&	3.73E-07	&	7.62E-05	\\\hline
				14	&	1.20E-03	&	1.50E-03	&	1.85E-04	&	1.50E-04	\\\hline
				24	&	2.37E-02	&	2.63E-02	&	6.28E-04	&	7.50E-04	\\\hline
				IEEE30	&	1.20E-03	&	5.03E-04	&	1.70E-04	&	1.67E-04	\\\hline
				30	&	5.15E-03	&	9.89E-03	&	1.57E-04	&	3.19E-04	\\\hline
				39	&	3.44E-04	&	1.63E-03	&	3.32E-04	&	1.64E-04	\\\hline
				57	&	9.58E-03	&	3.48E-02	&	2.06E-04	&	4.00E-04	\\\hline
				89	&	1.15E-04	&	3.68E-05	&	3.75E-05	&	6.60E-03	\\\hline
				118	&	2.31E-02	&	1.03E-02	&	8.32E-04	&	3.53E-03	\\\hline
				200	&	7.51E-03	&	8.68E-04	&	1.08E-04	&	2.01E-03	\\\hline
				300	&	1.66E-03	&	1.79E-03	&	1.34E-04	&	3.75E-04	\\\hline
				1354	&	3.44E-05	&	2.52E-05	&	1.61E-04	&	5.48E-03	\\\hline
				1888	&	2.76E-03	&	9.02E-03	&	1.53E-03	&	2.22E-02	\\\hline
				1951	&	1.65E-05	&	6.84E-06	&	1.34E-03	&	5.98E-03	\\\hline
				2383	&	6.63E-01	&	4.78E-01	&	4.38E-03	&	3.04E-03	\\\hline
				2736	&	1.22E-02	&	2.04E-02	&	2.30E-04	&	1.35E-04	\\\hline
				2737	&	4.12E-03	&	7.98E-03	&	2.08E-04	&	8.21E-05	\\\hline
				2746wop	&	8.11E-03	&	6.92E-03	&	1.74E-04	&	4.77E-04	\\\hline
				2746wp	&	5.77E-03	&	1.09E-02	&	1.67E-04	&	8.71E-05	\\\hline
				2848	&	7.72E-05	&	1.15E-04	&	8.20E-04	&	2.27E-02	\\\hline
				2868	&	2.92E-04	&	2.84E-04	&	1.03E-03	&	7.39E-03	\\\hline
				2869	&	3.43E-05	&	2.29E-05	&	1.10E-04	&	3.58E-03	\\\hline
				3012	&	9.25E-03	&	6.92E-03	&	1.41E-04	&	2.46E-04	\\\hline
				3120	&	2.96E-02	&	1.83E-02	&	3.10E-04	&	3.59E-04	\\\hline
				3375	&	7.36E-03	&	6.09E-03	&	7.96E-05	&	1.43E-04	\\\hline \hline
		\end{tabular}}
		\label{Table_LMP}
	\end{table}
	Detailed LMP and Q-LMP profiles for the IEEE 57-bus test system obtained from KNITRO and Algorithm~\ref{SLP_Meshed} are shown in Fig. \ref{fig_MP57}.

\section{Conclusion}\label{sec_conclusion}

This paper presented an SLP approach that exploits the structure of the alternative-form OPF by carefully constructing iteratively refined polyhedral outer approximations and supporting hyperplanes. Rigorous numerical evaluation on a wide range of test cases with up to 3375 buses shows that the method consistently converges to feasible high-quality solutions with robustness against the choice of starting point. While not claiming it is computationally superior to existing NLP solvers, the proposed method represents a competitive LP-based alternative, especially in settings where only LP solvers can be used or are practically deployed, e.g., market applications with shadow pricing of active and reactive power as well as congestion and losses. Another important benefit of the proposed SLP approach over NLP-based approaches is its potential straightforward extension to applications involving discrete decision variables, whereby powerful state-of-the-art MILP solvers can be used.


\appendices
\numberwithin{equation}{section}

\section{Branch current limit constraints}\label{sec_currentlimit}
	
	In the polar-form OPF, the branch current magnitude is expressed as
	\begin{subequations}\label{eq_currentmag}
		\begin{align}
		\hspace{-0.22cm}	i_{ij} = \sqrt{ a_{1}v_{i}^{2} + a_{2}v_{j}^{2} - 2v_{i}v_{j}\left( a_{3}\cos{\theta_{ij}} - a_{4}\sin{\theta_{ij}} \right)  },\\
		\hspace{-0.22cm}	i_{ji} = \sqrt{ a_{2}v_{i}^{2} + a_{1}v_{j}^{2} - 2v_{i}v_{j}\left( a_{3}\cos{\theta_{ij}} + a_{4}\sin{\theta_{ij}} \right)  },
		\end{align}
	\end{subequations}
	where $a_{1} = \left( \tilde{g}^{\rm c^*}_{ij} \right)^{2} $, $a_{2} = \left( \tilde{g}^{*}_{ij} \right)^{2} $, $a_{3} = \Re\{ \tilde{g}^{\rm c^*}_{ij} \}\Re\{ \tilde{g}^{*}_{ij} \} + \Im\{ \tilde{g}^{\rm c^*}_{ij} \}\Im\{ \tilde{g}^{*}_{ij} \} $, and $a_{4} = \Re\{ \tilde{g}^{\rm c^*}_{ij} \}\Im\{ \tilde{g}^{*}_{ij} \} - \Im\{ \tilde{g}^{\rm c^*}_{ij} \}\Re\{ \tilde{g}^{*}_{ij} \} $. The limit constraint on the branch current magnitude can now be written as
	\begin{align}
		i_{ij} & \le \overline{I}_{ij}, & ij \in \mathcal{L} \cup \mathcal{L}_{\rm t},
	\end{align}
	where $\overline{I}_{ij} = \overline{I}_{ji}$ (pu) is the current rating of branch $ij$.
	
	Constraints \eqref{eq_currentmag} are nonconvex and, as shown in \cite{Shchetinin2019_CurrentLimits}, linearizing them is a daunting task. However, in the \emph{alternative-form} OPF, those constraints can be replaced by \emph{linear} equivalents. In more detail, by defining $\ell_{ij} := i_{ij}^{2}$, and after substituting \eqref{opf_definingW} in \eqref{eq_currentmag} and squaring both sides of the resulting equation, the branch current magnitude and associated limit become
	\begin{subequations}\label{eq_currentmag_alt}
		\begin{align}
			&\ell_{ij} = a_{1}w_{i} + a_{2}w_{j} - 2a_{3}w_{ij}^{\rm r} + 2a_{4}w_{ij}^{\rm i} , & ij \in \mathcal{L}&, \\
			&\ell_{ji} = a_{2}w_{i} + a_{1}w_{j} - 2a_{3}w_{ij}^{\rm r} - 2a_{4}w_{ij}^{\rm i}, & ij \in \mathcal{L}&, \\
			&\ell_{ij} \le \overline{I}_{ij}^{2}. & \hspace{-1cm} ij \in \mathcal{L} \cup \mathcal{L}&_{\rm t}
		\end{align}
	\end{subequations}
	The original current magnitude $i_{ij}$ can be straightforwardly obtained from the solution by taking the square root of $\ell_{ij}$ in a computationally cheap post-processing step.

\section{Proof of \cref{opf_socr_slp_Convergence}}\label{opf_socr_slp_Proof}
	
	\begin{proof}[\unskip\nopunct]
		If $x^{(k)}$ minimizes $f_{0}(x)$ on $\mathcal{S}^{(k)}$ then it must satisfy the inequalities
		\begin{align}\label{convergence_proof}
			f_{i}(x^{(\kappa)}) + \nabla f_{i}(x^{(k)})^{T}(x^{(k)} - x^{(\kappa)}) \leq 0, \nonumber \\
			\kappa=0,\ldots,k-1, \ i=1,\ldots,m.
		\end{align}
		Additionally, if $x^{(k)}$ contains a subsequence that converges to a point in $\mathcal{G}$ then every $f_{i}(x^{(k)}), i=1,\ldots,m\ $ must contain a subsequence that converges to a value that is less than or equal to zero. If there exists a finite constant $K$ such that $\max \left( \left( \left\| \nabla f_{i}(x^{(k)}) \right\| \right)_{i=1,\ldots,m} \right) \le K $ for all $x \in \mathcal{S}^{(0)}$, and an $r > 0$ independent of $k$, such that
		\begin{align*}
			\resizebox{\linewidth}{!}{$r \le f_{i}\left(x^{(\kappa)}\right) \le \nabla f_{i}(x^{(k)})^{T}(x^{(\kappa)} - x^{(k)})\leq K\left\| x^{(\kappa)} - x^{(k)} \right\|,$}\\
			\kappa=0,\ldots,k-1, \ i=1,\ldots,m,
		\end{align*}
		then one can find a subsequence of iteration numbers where
		\begin{align*}
			\left\| x^{(u)} - x^{(v)} \right\| \le r/K, \ \forall u < v,
		\end{align*}
		so that $x^{(k)}$ does \emph{not} contain a Cauchy sequence. However, because $\mathcal{S}^{(0)}$ is compact, no such $r$ exists, and $x^{(k)}$ therefore contains a subsequence that converges to a point $x^{\dagger} \in \mathcal{S}^{(0)}$. Finally, it follows from \eqref{convergence_proof} that the corresponding subsequence of every $f_{i}(x^{(k)}), i=1,\ldots,m\ $ converges to a value that is less than or equal to zero, which entails that $x^{\dagger} \in \mathcal{G}$ and $x^{\star}=x^{\dagger}$.
	\end{proof}

\section{SLP example}\label{slp_example}

\begin{figure}[t!]
		\centering
		\begin{subfigure}{\columnwidth}
			\includegraphics[width=\columnwidth] {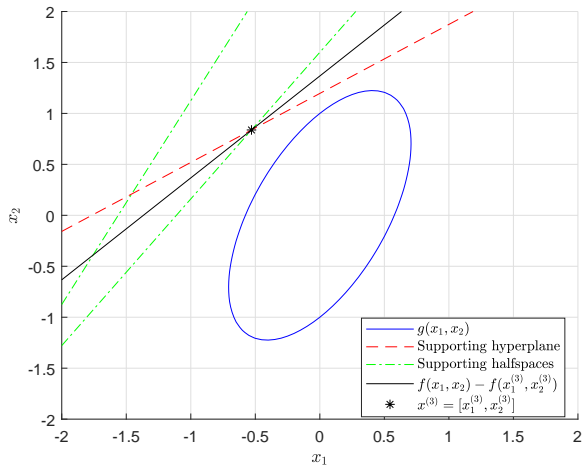}
			\caption{$k=3$.}
			\label{fig_SLPk3}
		\end{subfigure}
		\begin{subfigure}{\columnwidth}
			\includegraphics[width=\columnwidth] {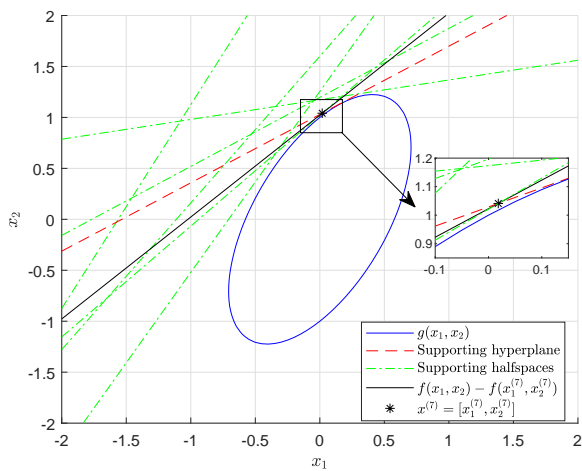}
			\caption{$k=7$.}
			\label{fig_SLPk7}
		\end{subfigure}
		\caption{Geometric interpretation of the SLP algorithm at $k=3$ and $k=7$.}
		\label{fig_SLPexample}
	\end{figure}

\begin{figure}[t!]
		\centering{
			\includegraphics[width=1.0\columnwidth] {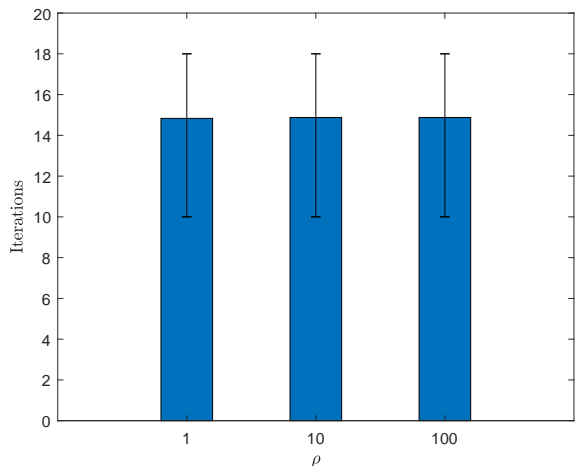}}
		\caption{Average number of iterations to convergence of the proposed SLP algorithm initialized from 1000 different starting points for $\rho = 1,10,100$. The error bars show the maximum and minimum number of iterations to convergence.}
		\label{fig_SLProbustness}
	\end{figure}

	Consider the nonconvex problem
	\begin{subequations}\label{example_nlp}
		\begin{align}
			\underset {\substack{x}} 
			{\mbox{ minimize}} \quad & f\left(x\right) = x_{1} - x_{2} & \label{example_nlp_objective}\\
			\text{ subject to} \quad g\left(x\right) & = 3x_{1}^{2} - 2x_{1}x_{2} + x_{2}^{2} -1 = 0, & \label{example_nlp_constraint}\\
			-2 \le & x_{1},x_{2} \le 2, & \label{example_nlp_box}
		\end{align}
	\end{subequations}
	whose optimal solution is $x^{*}=[x_{1}^{*},x_{2}^{*}]^{T}=[0,1]^{T}$. The novel SLP algorithm in this work solves the problem
	\begin{subequations}\label{example_slp}
		\begin{align}
			\underset {\substack{x}} 
			{\mbox{ minimize}} \quad & f\left(x\right) + \rho r & \label{example_slp_objective}\\
			\hspace{-0.5cm} \text{ subject to} \ g\left(x^{(k)}\right) + & \nabla g\left(x^{(k)}\right)^{T} \left( x-x^{(k)}\right) = - r, \hspace{-1cm} & \label{example_slp_hyperp}\\
			g\left(x^{(\kappa)}\right) + & \nabla g\left(x^{(\kappa)}\right)^{T} \left( x-x^{(\kappa)}\right) \le 0, \hspace{-1cm} & \nonumber \\
			 & \hspace{1.75cm} \kappa=0,\ldots,k-1 & \label{example_slp_halfsp}\\
			-2 \le & x_{1},x_{2} \le 2, & \label{example_slp_box} \\
			r & \ge 0, & \label{example_slp_slack}
		\end{align}
	\end{subequations}
	at each iteration $k$ until convergence, measured by $ \left| g\left(x^{(k)}\right) \right| \le \epsilon$. A geometric interpretation of the SLP algorithm with $\rho=10$ and from starting point $x^{(0)}=[-2,2]^{T}$ is shown in Fig.~\ref{fig_SLPk3} and Fig.~\ref{fig_SLPk7} at $k=3$ and $k=7$, respectively.
	In Fig.~\ref{fig_SLPexample}, supporting hyperplanes \eqref{example_slp_hyperp} are shown in red and supporting halfspaces \eqref{example_slp_halfsp} are shown in green. The contour lines of $f\left( x \right) - f \left( x^{(k)} \right)$ and the nonconvex constraint $g\left(x\right) = 0$ are shown in black and blue, respectively. The simultaneous inclusion of the supporting hyperplanes \eqref{example_slp_hyperp} and the supporting halfspaces \eqref{example_slp_halfsp} is essential for convergence. Without the hyperplanes in \eqref{example_slp_hyperp} the algorithm converges to the solution of a convex relaxation of Problem~\ref{example_nlp}, \emph{i.e.,} one with $g\left(x\right) \le 0$ instead of $g\left(x\right) = 0$. Moreover, without the halfspaces in \eqref{example_slp_halfsp} the algorithm fails to converge altogether. Furthermore, without the positive slack variable $r$ and the second term in \eqref{example_slp_objective} some starting points, such as $x^{(0)}=[0,0]^{T}$ in this case, can lead to an infeasible LP Problem~\ref{example_slp}. The robustness against the choice of starting point bestowed by $\rho r$ and \eqref{example_slp_slack} is shown in Fig. \ref{fig_SLProbustness} for $\epsilon = 10^{-4}$, three choices of $\rho$, and 1000 starting points drawn from a uniform distribution in the interval $(-10,10)$ for both $x_{1}^{(0)}$ and $x_{2}^{(0)}$, \emph{i.e.,} $x_{1}^{(0)},x_{2}^{(0)} \sim \mathcal{U}\left( -10,10\right) $.
	The error plot in Fig. \ref{fig_SLProbustness} demonstrates a small variation in the number of iterations to convergence for $\rho = 1,10,100$ over all the 1000 starting points, thus corroborating the robustness of the proposed SLP algorithm on Problem~\ref{example_nlp}.
	
\begin{figure}[t!]
		\centering{
			\includegraphics[width=1.0\columnwidth] {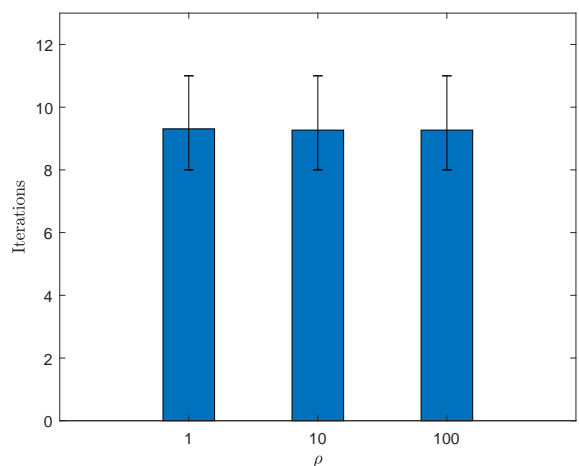}}
		\caption{Average number of iterations to convergence of the proposed SLP algorithm with the projection of point $x^{(k)}$ on the set defined by $g\left(x\right)=0$. The algorithm is again initialized from 1000 different starting points for $\rho = 1,10,100$. The error bars show the maximum and minimum number of iterations to convergence.}
		\label{fig_SLProbustness_proj}
	\end{figure}

	More interestingly, the SLP algorithm applied to Problem~\ref{example_nlp} converges much faster if the projection of point $x^{(k)}$ on the set defined by $g\left(x\right)=0$ is used instead of the first-order Taylor series approximation of $g\left(x\right)$ directly at $x^{(k)}$ (See Section \ref{subsec_SLPforSOCR} for more details on the projection). 
	\begin{figure}[t!]
		\centering{
			\includegraphics[width=1.0\columnwidth] {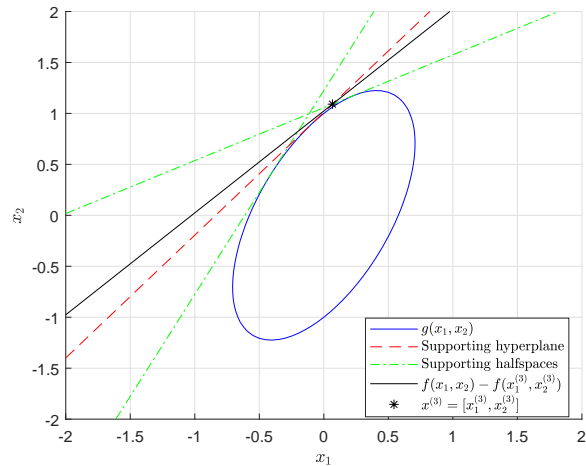}}
		\caption{Geometric interpretation of the SLP algorithm with the projection of point $x^{(k)}$ on the set defined by $g\left(x\right)=0$ at $k=3$.}
		\label{fig_SLPk3_proj}
	\end{figure}
	Fig.~\ref{fig_SLPk3_proj} shows that the modified SLP algorithm is already close to the optimal solution after only 3 iterations. In fact, as shown in Fig. \ref{fig_SLProbustness_proj}, the proposed SLP algorithm with the projection of point $x^{(k)}$ on the set defined by $g\left(x\right)=0$ is on average $40\%$ faster compared to directly applying the first-order Taylor series approximation of $g\left(x\right)$ at $x^{(k)}$ (9 iterations on average for the modified SLP algorithm compared to 15 for the one with the first-order Taylor series approximation of $g\left(x\right)$ directly at $x^{(k)}$).
	Finally, it is also evident form Fig.~\ref{fig_SLPexample} and Fig.~\ref{fig_SLPk3_proj} that not all the supporting halfspaces are necessary for convergence. Some of the older inactive halfspaces can be removed for an improvement in computational speed.
	
\bibliographystyle{IEEEtran}
{\footnotesize
	\bibliography{SLPforOPF}}

	\begin{IEEEbiography}[{\includegraphics[width=1in,height=1.25in,clip,keepaspectratio]{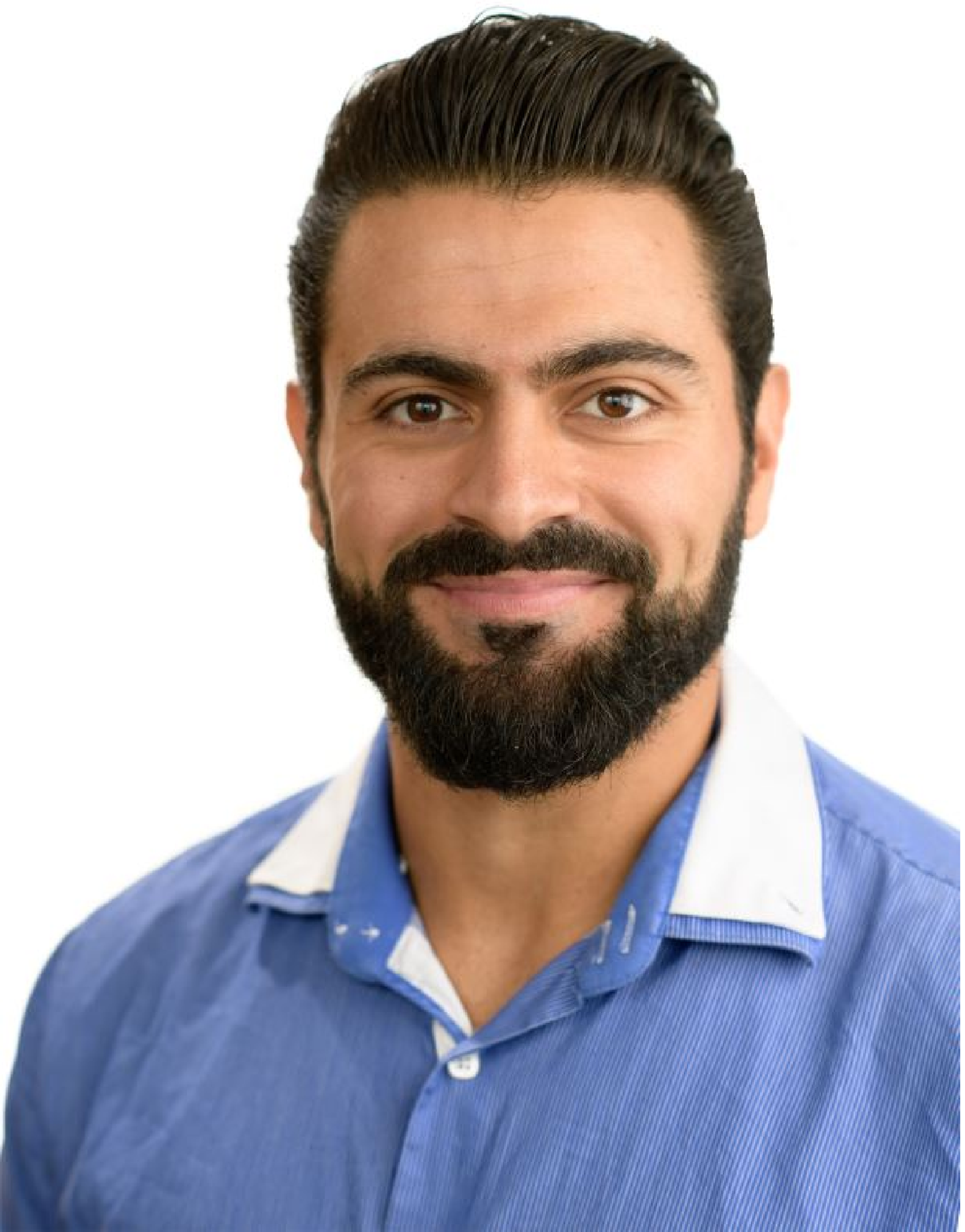}}]{Sleiman Mhanna}
		(S'13--M'16) received the B.Eng. degree (with high distinction) from the Notre Dame University, Lebanon, and the M.Eng. degree from the American University of Beirut, Lebanon, in 2010 and 2012, respectively, both in electrical engineering. He received the Ph.D. degree from the School of Electrical and Information Engineering, Centre for Future Energy Networks, University of Sydney, Australia, in 2016. He is currently a Research Fellow at the Department of Electrical and Electronic Engineering, The University of Melbourne, Australia. His research interests include computational methods for integrated multi-energy systems, decomposition methods, and demand response.
	\end{IEEEbiography}
	
	\begin{IEEEbiography}[{\includegraphics[width=1in,height=1.25in,clip,keepaspectratio]{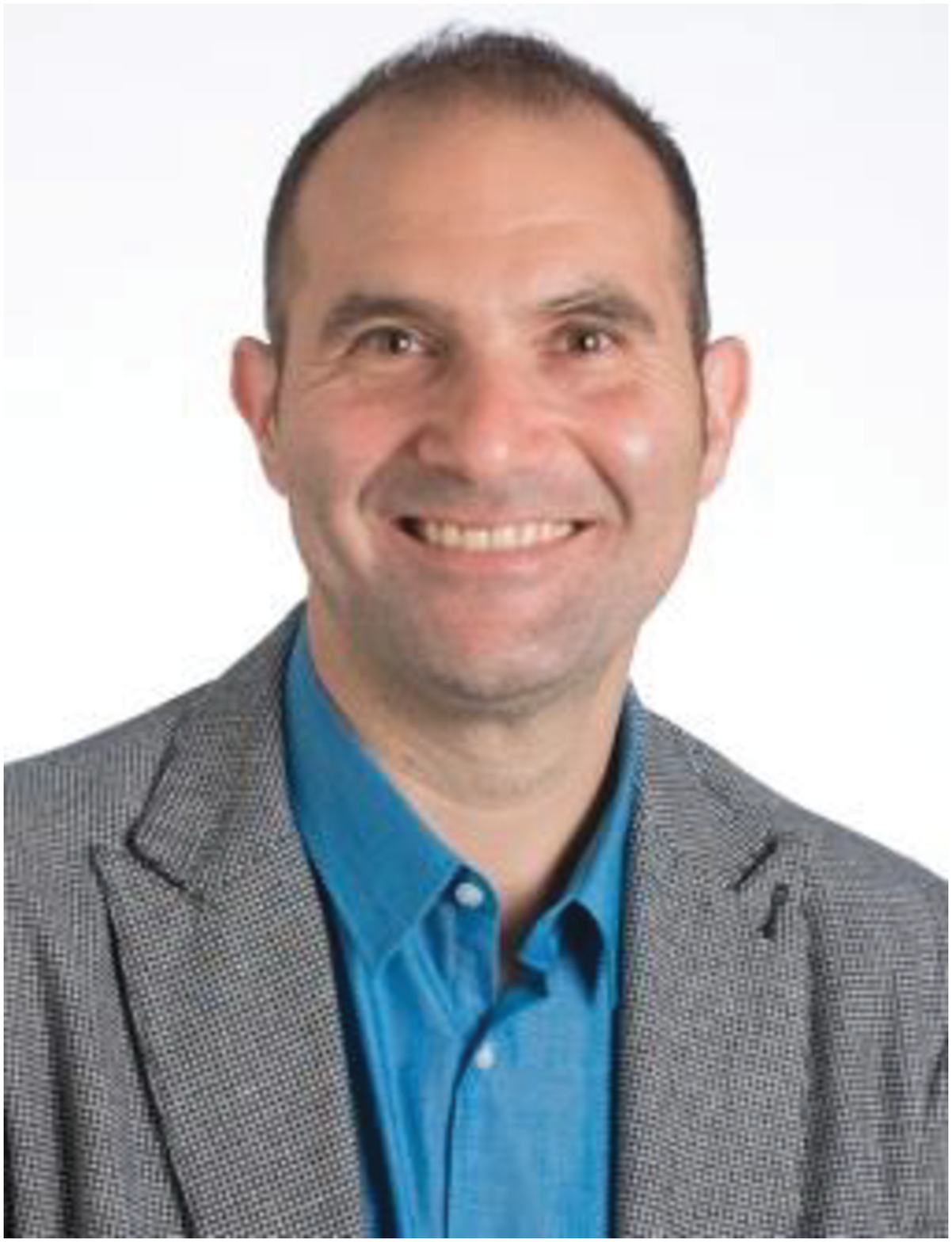}}]{Pierluigi Mancarella}
		(SM) is Chair Professor of Electrical Power Systems at The University of Melbourne, Melbourne, Australia, and Professor of Smart Energy Systems at The University of Manchester, Manchester, UK. His research interests include techno-economic modeling of integrated multi-energy systems; security, reliability and resilience of future networks; and energy infrastructure planning under uncertainty. Pierluigi is an Editor of the IEEE Transactions on Power Systems, IEEE Transactions on Smart Grid, and IEEE Systems Journal, and an IEEE Power and Energy Society Distinguished Lecturer.
	\end{IEEEbiography}

\end{document}